\documentclass[12pt]{imsart}
\usepackage{amsthm,amsmath}
\usepackage[colorlinks,citecolor=blue,urlcolor=blue]{hyperref}

\usepackage{amsmath,amsfonts,amssymb,amsthm}

\usepackage{enumitem}

\RequirePackage{amsthm,amsmath,amsfonts,amssymb}
\RequirePackage[numbers]{natbib}
\RequirePackage[colorlinks,citecolor=blue,urlcolor=blue]{hyperref}
\RequirePackage{graphicx}
\usepackage{lineno,hyperref}

\startlocaldefs
\modulolinenumbers[5]

\usepackage[utf8]{inputenc}
\usepackage[english]{babel}

\usepackage[margin=1in]{geometry}
\usepackage{bm}
\usepackage{graphicx}
\usepackage{amssymb,amsmath,amsthm,mathrsfs}
\usepackage{epstopdf}
\usepackage{algorithm}
\usepackage{amsfonts,delarray}
\usepackage{algorithm,algcompatible,amsmath}
\usepackage{rotating,color}
\usepackage{subfigure}
\usepackage{multirow}
\usepackage{titlesec,textcase}
\usepackage{longtable}
\usepackage{bbm}
\usepackage{rotating}

\newtheorem{Theorem}{Theorem}[section]
\newtheorem{Lemma}[Theorem]{Lemma}

\newtheorem{Proposition}[Theorem]{Proposition}

\theoremstyle{definition}

\RequirePackage[normalem]{ulem} 

\definecolor{rp}{RGB}{83,54,106}

\def\boxit#1{\vbox{\hrule\hbox{\vrule\kern6pt\vbox{\kern6pt#1\kern6pt}\kern6pt\vrule}\hrule}}

\endlocaldefs

\begin{document}
\begin{frontmatter}
\title{On the R\'{e}nyi index of random graphs}

\runtitle{heterogeneity of network }
\runauthor{ }
\begin{aug}

\author[A]{\fnms{Mingao} \snm{Yuan}\ead[label=e1]{mingao.yuan@ndsu.edu}}



\address[A]{Department of Statistics,
North Dakota State University,
\printead{e1}}
\end{aug}

\begin{abstract}
Networks (graphs) permeate scientific fields such as biology, social science, economics, etc. Empirical studies have shown that real-world networks are often heterogeneous, that is, the degrees of nodes do not concentrate on a number.  Recently,
the R\'{e}nyi index was tentatively used to measure network heterogeneity.  However, the validity of the R\'{e}nyi index in network settings is not theoretically justified. 
In this paper, we study this problem. We 
derive the limit of the R\'{e}nyi index of a heterogeneous Erd\"{o}s-R\'{e}nyi random graph and a power-law random graph, as well as the convergence rates. Our results show that the  Erd\"{o}s-R\'{e}nyi random graph has asymptotic R\'{e}nyi index zero and the power-law random graph (highly heterogeneous) has asymptotic  R\'{e}nyi index one. In addition, the limit of the R\'{e}nyi index increases as the graph gets more heterogeneous. These results theoretically justify the R\'{e}nyi index is a reasonable statistical measure of network heterogeneity. We also evaluate the finite-sample performance of  the R\'{e}nyi index by simulation.
\end{abstract}



\end{frontmatter}

\section{Introduction}
\label{S:1}

A network (graph) consists of a set of individuals (nodes) and a set of interactions (edges) between individuals. It has been widely used to model and analyze many complex systems. For example, in social science and economics, networks play a central role in the transmission of information, the trade of many goods and services, and determining how diseases spread, etc. \cite{C18,KRTB22,REE08}; in biology, network is a method of representing  the physical contacts between proteins \cite{CY06}. In the past decade, network data analysis has been a primary research topic in statistics and machine learning \cite{A18,ACB13,BS16,GZFA,YS21,YS21b}.

In many fields of science and engineering, one of the elemental problems is to measure the statistical heterogeneity of datasets. For instance, in statistical physics the entropy was devised to measure the randomness of systems (\cite{K07}). In economics, various inequality indices were designed to gauge the evenness of the distribution of wealth in human populations(\cite{C89}). Motivated by entropy and inequality indices, \cite{E11} recently introduced the R\'{e}nyi index to measure statistical heterogeneity
of probability distributions defined on the positive half-line. The R\'{e}nyi index takes values in the range $[0,1]$. Larger value represents higher level of heterogeneity. Its properties were systematically studied in \cite{E11,ES12} and the R\'{e}nyi index of several well-known distributions (such as Pareto distribution, Gamma distribution, Beta distribution, etc.) are calculated in \cite{E11,ES12}.

Empirical studies have shown that many real-world networks are heterogeneous, that is, the degrees of individuals do not concentrate on a number \cite{CSN09,N03,VHHK19}. It is important to be able to compare networks according to heterogeneity that they exhibit, and thus to
have a stable summary statistic that provides
insight into the structure of a network.
Recently the R\'{e}nyi index was tentatively used to measure heterogeneity of financial networks and interesting findings were obtained \cite{ NSL16,NS19,NS21,N21}.
However, the validity of the R\'{e}nyi index in network  settings is not theoretically justified, and some of the fundamental questions are not studied in \cite{NSL16}. For instance, whether the R\'{e}nyi index of a homogeneous network is actually close to zero, whether the R\'{e}nyi index of heterogeneous network is indeed large, and how the R\'{e}nyi index depends on network model parameters. 

In this paper, we shall answer the above mentioned questions and provide a theoretical justification for the R\'{e}nyi index as a network heterogeneity measure. To this end, we derive the limit of the R\'{e}nyi index of a heterogeneous Erd\"{o}s-R\'{e}nyi random graph and a power-law random graph, as well as the convergence rates. Based on our results, the Erd\"{o}s-R\'{e}nyi random graph (homogeneous) has asymptotic R\'{e}nyi index zero, while the well-known  power-law random graph (highly heterogeneous) has asymptotic   R\'{e}nyi index one. Moreover, the limit of the R\'{e}nyi index explicitly depends on model parameters, from which it is clear that the R\'{e}nyi index increases as the model gets more heterogeneous.
These results theoretically justify the R\'{e}nyi index is a reasonable statistical measure of network heterogeneity. In addition, we run simulations to evaluate finite-sample performance of the R\'{e}nyi index.

The structure of the article is as follows. In Section \ref{main} we collect the main results. Specifically, in Section \ref{mainsub}, we present the limit of the R\'{e}nyi index of a heterogeneous Erd\"{o}s-R\'{e}nyi random graph; in Section \ref{mainpower},  we present the limit of the R\'{e}nyi index of a power-law random graph. Simulation studies are given in Section \ref{simu} . All proofs are deferred to Section \ref{proof}.

\medskip

{\bf Notation}:  Let $c_1,c_2$ be two positive constants. For two positive  sequences $a_n$, $b_n$, denote $a_n\asymp b_n$ if $c_1\leq \frac{a_n}{b_n}\leq c_2$; denote $a_n=O(b_n)$ if $\frac{a_n}{b_n}\leq c_2$; $a_n=o(b_n)$ if $\lim_{n\rightarrow\infty}\frac{a_n}{b_n}=0$. Let $X_n$ be a sequence of random variables. We use $X_n\Rightarrow F$ to denote $X_n$ converges in distribution to a probability distribution $F$. $X_n=O_P(a_n)$ means $\frac{X_n}{a_n}$ is bounded in probability.  $X_n=o_P(a_n)$ means $\frac{X_n}{a_n}$ converges to zero in probability. $\mathcal{N}(0,1)$ stands for the standard normal distribution. Let $I[E]$ be the indicator function of an event $E$. We adopt the convention $0 \log 0 = 0$. Let $n$ be a positive integer.

\section{The R\'{e}nyi index of network}\label{main}

A graph or network 
 $\mathcal{G}$ consists of a pair $(\mathcal{V},\mathcal{E})$, where $\mathcal{V}=[n]:=\{1,2,\dots,n\}$ denotes the set of vertices and $\mathcal{E}$ denotes the set of edges. 
 For $i<j$, denote $A_{ij}=1$ if $\{i,j\}\in\mathcal{E}$ is an edge and $A_{ij}=0$ otherwise. Suppose  $A_{ii}=0$, that is,  self loops are not allowed. Then the symmetric matrix $A=(A_{ij})\in\{0,1\}^{\otimes n^2}$ is called the adjacency matrix of graph $\mathcal{G}$. A graph is said to be random if the elements $A_{ij} (1\leq i<j\leq n)$ are random.

Given a positive constant $\alpha$, the R\'{e}nyi index of a graph (\cite{E11,ES12,NSL16,NS19}) is defined as 
\begin{equation}\label{renyiindex}
  \mathcal{R}_{\alpha}=
    \begin{cases}
     1-\left[\frac{1}{n}\sum_{i=1}^n\left(\frac{d_i}{d}\right)^{\alpha}\right]^{\frac{1}{1-\alpha}},  & \text{if $\alpha\neq 1$;}\\
      1-\exp\left(-\frac{1}{n}\sum_{i=1}^n\frac{d_i}{d}\log\frac{d_i}{d}\right), & \text{if $\alpha=1$,}
    \end{cases}       
\end{equation}
where $d_i$ is the degree of node $i$, that is, $d_i=\sum_{j\neq i}A_{ij}$ and $d$ is the average of degree, that is, $d=\frac{\sum_{i=1}^nd_i}{n}$. The R\'{e}nyi index includes several popular indexes as a special case.  When $\alpha=1$, the R\'{e}nyi index $\mathcal{R}_{1}$ is a function of the Theil's index. When $\alpha=2$, the R\'{e}nyi index $\mathcal{R}_{2}$ is function of the Simpson's index. For $0<\alpha\leq 1$ the R\'{e}nyi index $\mathcal{R}_{\alpha}$ is the Atkinson's index. The parameter $\alpha$ allows researchers to tune the R\'{e}nyi index to be more sensitive to different populations. In practice, commonly used values are $\alpha=1,2,3$ (\cite{NSL16,NS19}).

\begin{Proposition}\label{prop}
For any fixed $\alpha>0$,
 the R\'{e}nyi index $ \mathcal{R}_{\alpha}$ is between 0 and 1.
\end{Proposition}

The R\'{e}nyi index takes values in $[0,1]$. It is tentatively used to measure degree heterogeneity of graphs (\cite{NSL16}). We shall derive an asymptotic expression  of the R\'{e}nyi index of two random graphs. 
 Note that $\mathcal{R}_{\alpha}$ is a non-linear function of the degrees $d_i\ (1\leq i\leq n)$, the degrees are not independent and may not be identically distributed. This fact make studying asymptotic properties of the R\'{e}nyi index a non-trivial task.

\subsection{The R\'{e}nyi index of a heterogeneous Erd\"{o}s-R\'{e}nyi random graph}\label{mainsub}

In this section, we study the asymptotic  R\'{e}nyi index of a heterogeneous Erd\"{o}s-R\'{e}nyi random graph.  Let $f(x,y)$ be a symmetric function from $[0,1]^2$ to $[0,1]$.
Define the heterogeneous Erd\"{o}s-R\'{e}nyi random graph $\mathcal{G}(n,p_n, f)$ as
\[\mathbb{P}(A_{ij}=1)=p_n f\left(\frac{i}{n},\frac{j}{n}\right),\]
where $p_n\in[0,1]$ may depend on $n$ and $A_{ij}\ (1\leq i<j\leq n)$ are independent. If $f\equiv c$ for some constant $c$, then $\mathcal{G}(n,p_n, f)$ is simply the Erd\"{o}s-R\'{e}nyi random graph with edge appearance probability $cp_n$. For non-constant $f$, the random graph $\mathcal{G}(n,p_n, f)$ is a heterogeneous version of the Erd\"{o}s-R\'{e}nyi graph. 
The spectral properties of this random graph have been extensively studied in \cite{CHHS21,CHHS20,CCH20}. We point out that our proof strategy works for other graph models such as the $\beta$-model in \cite{RPF13} and the degree-corrected model in \cite{KN11} with mild modifications.

\subsubsection{Asymptotic R\'{e}nyi index when $\alpha\neq 1$}

In this subsection, we study asymptotic R\'{e}nyi index of $\mathcal{G}(n,p_n, f)$ with $\alpha\neq 1$.
For convenience, denote
\[f_{ij}=f\left(\frac{i}{n},\frac{j}{n}\right),\hskip 1cm f_i=\frac{1}{n}\sum_{j\neq i}^nf\left(\frac{i}{n},\frac{j}{n}\right),\ \ \ \ \lambda_{k,l}=\frac{\sum_{i\neq j}f_i^{k}f_{ij}^l}{n^2}.\]
 Note that $f_{ij}$, $f_i$ and $\lambda_{k;l}$ depend on $n$. We will focus on $f(x,y)\geq \epsilon$ for a constant $\epsilon\in(0,1)$ as assumed in \cite{CHHS20}. Later we will provide examples of such functions.

\begin{Theorem}\label{theorem:1} Let $\alpha\neq 1$ be a fixed positive constant, $np_n\rightarrow\infty$ and $f(x,y)\geq \epsilon$ for some constant $\epsilon\in(0,1)$. Then the R\'{e}nyi index $\mathcal{R}_{\alpha}$ of $\mathcal{G}(n,p_n,f)$ has the following expression
\begin{equation}\label{renyiasymp}
   \mathcal{R}_{\alpha}=1-\left[\frac{\lambda_{\alpha,0}}{\left(\lambda_{0,1}+O_P\left(\frac{1}{n\sqrt{p_n}}\right)\right)^{\alpha}}+O_P\left(\frac{1}{np_n}\right)\right]^{\frac{1}{1-\alpha}},
\end{equation}
and the error rates $\frac{1}{np_n}$ and $\frac{1}{n\sqrt{p_n}}$ cannot be improved. Asymptotically, $\mathcal{R}_{\alpha}$ has the following concise  expression:
\begin{equation}\label{renyiind}
\mathcal{R}_{\alpha}=1-\left(\frac{\lambda_{\alpha,0}}{\lambda_{0,1}^{\alpha}}\right)^{\frac{1}{1-\alpha}}+o_P(1).
\end{equation}
\end{Theorem}
\medskip

Theorem \ref{theorem:1} provides   an asymptotic expression of $\mathcal{R}_{\alpha}$ as an explicit function of $\alpha$ and the model parameter $f$, along with the error rates.
It is interesting that $\mathcal{R}_{\alpha}$ mainly depends on $f$ and $\alpha$ through the ratio $\frac{\lambda_{\alpha,0}}{\lambda_{0,1}^{\alpha}}$. The quantities $\lambda_{\alpha,0}$ and $\lambda_{0,1}^{\alpha}$ may or may not converge to some limits as $n$ goes to infinity. Later we will present two examples where $\lambda_{\alpha,0}$ and $\lambda_{0,1}^{\alpha}$ converge.

We point out that even though empirical degree distributions are widely studied in literature, it is not immediately clear how to obtain the asymptotic expression of the R\'{e}nyi index as in Theorem \ref{theorem:1} from the empirical degree distributions. Specifically, suppose $Y_n$ is a random variable with distribution $F_{emp}(x)=\frac{1}{n}\sum_{i=1}^nI[d_i\leq x]$. The term $\frac{1}{n}\sum_{i=1}^nd_i^{\alpha}$ in the R\'{e}nyi index is equal to $\mathbb{E}(Y_n^{\alpha}|d_1,\dots,d_n)$. Suppose $F_{emp}(x)$ converges almost surely or in probability to some distribution function $F(x)$ and let $Y$ follow the distribution $F(x)$. The convergence of $F_{emp}(x)$ to $F(x)$ does not necessarily imply the convergence of $\mathbb{E}(Y_n^{\alpha}|d_1,\dots,d_n)$ to $\mathbb{E}(Y^{\alpha})$ for arbitrary $\alpha>0$. Generally speaking, uniform integrability conditions are required to guarantee the convergence of $\mathbb{E}(Y_n^{\alpha}|d_1,\dots,d_n)$ to $\mathbb{E}(Y^{\alpha})$. Note that $\mathbb{E}(Y_n^{\alpha}|d_1,\dots,d_n)$ is random. It is not immediately clear what kind of uniform integrability conditions are needed. Moreover, even if we can conclude that $\mathbb{E}(Y_n^{\alpha}|d_1,\dots,d_n)$ converges to $\mathbb{E}(Y^{\alpha})$ by assuming some uniform integrability conditions,  it does not provide the error rates (that cannot be improved) as in Theorem \ref{theorem:1}.

Next we provide two examples of random graphs satisfying the conditions of Theorem \ref{theorem:1} and  calculate the ratio explicitly.

The first example is the Erd\"{o}s-R\'{e}nyi random graph, that is, $f(x,y)\equiv 1$. Since each node of the Erd\"{o}s-R\'{e}nyi random graph has the same average degree,  the Erd\"{o}s-R\'{e}nyi graph is homogeneous.
It is clear that $\frac{\lambda_{\alpha,0}}{\lambda_{0,1}^{\alpha}}=1+o(1)$, hence $\mathcal{R}_{\alpha}=o_P(1)$. This shows the R\'{e}nyi index of homogeneous network is actually close to zero.

Now we provide a family of non-constant $f(x,y)$ that is bounded away from zero. This model can attain any heterogeneity level, that is, the limit of $\frac{\lambda_{\alpha,0}}{\lambda_{0,1}^{\alpha}}$ can take any value in $(0,1)$.
Let $f(x,y)=e^{-\kappa x}e^{-\kappa y}$ with a non-negative constant $\kappa$. Then $e^{-2\kappa}\leq f(x,y)\leq 1$ for  $0\leq x,y\leq 1$. Intuitively, smaller $\kappa$ would produce less heterogeneous models. In the extreme case $\kappa=0$, the random graph is simply the Erd\"{o}s-R\'{e}nyi random graph. Given a function $f$, denote the expected degree of node $i$ as $\mu_i:=p_n\sum_{j\neq i}f_{ij}$ . Then for $f(x,y)=e^{-\kappa x}e^{-\kappa y}$, $\mu_i$ is equal to
\[np_ne^{-\kappa \frac{i}{n}}(1-e^{-\kappa}+o(1)).\]
Note that $\frac{\mu_{1}}{\mu_{n}}=e^{\kappa\left(1-\frac{1}{n}\right)}$. Large $\kappa$ will enlarge the difference between the degrees of node $1$ and node $n$. Hence, the random graph with larger $\kappa$ should be more heterogeneous. Simple calculations yield
\[\lambda_{\alpha,0}=\left(\frac{1}{\kappa}-\frac{1}{\kappa e^{\kappa}}\right)^{\alpha}\left(\frac{1}{\kappa\alpha}-\frac{1}{\kappa\alpha e^{\kappa\alpha}}\right)+o(1),\hskip 1cm \lambda_{0,1}=\left(\frac{1}{\kappa}-\frac{1}{\kappa e^{\kappa}}\right)^2+o(1).\]
Plugging them into (\ref{renyiind}) yields
\begin{equation}\label{lrenyi}
\mathcal{R}_{\alpha}=1-\left(\frac{(e^{\kappa\alpha}-1)\kappa^{\alpha-1}}{\alpha(e^{\kappa}-1)^{\alpha}}\right)^{\frac{1}{1-\alpha}}+o_P(1),\ \ \ \ \alpha>0,\ \alpha\neq 1.
\end{equation}
Note that $\lim_{\kappa\rightarrow\infty}\left(\frac{(e^{\kappa\alpha}-1)\kappa^{\alpha-1}}{\alpha(e^{\kappa}-1)^{\alpha}}\right)^{\frac{1}{1-\alpha}}=0$ and $\lim_{\kappa\rightarrow0^+}\left(\frac{(e^{\kappa\alpha}-1)\kappa^{\alpha-1}}{\alpha(e^{\kappa}-1)^{\alpha}}\right)^{\frac{1}{1-\alpha}}=1$ for any $ \alpha>0$ and $\alpha\neq 1$. Asymptotically,  $\mathcal{R}_{\alpha}$ with large $\kappa$ would be close to 1 and $\mathcal{R}_{\alpha}$ with small $\kappa$ would be close to 0. This justifies that the  R\'{e}nyi index of heterogeneous graph is actually non-zero.
In addition, the limit of $\mathcal{R}_{\alpha}$ can assume any value in $(0,1)$ by changing $\kappa$. In this sense, this random graph can achieve any heterogeneity level with suitably selected $\kappa$.

\subsubsection{Asymptotic R\'{e}nyi index when $\alpha=1$}

In this subsection, we study the asymptotic  R\'{e}nyi index of $\mathcal{G}(n,p_n, f)$ with $\alpha=1$. For convenience, denote
\[f_{ij}=f\left(\frac{i}{n},\frac{j}{n}\right),\hskip 3mm \mu_i:=p_n\sum_{j\neq i}f_{ij},\hskip 3mm l_i=\log\left(\frac{\mu_i}{np_n\lambda_{0,1}}\right), \hskip 3mm  s_k=\sum_{i<j}(2+l_i+l_j)^kf_{ij}(1-p_nf_{ij}), \]

\begin{Theorem}\label{theorem:2} Let $\mathcal{G}(n,p_n,f)$ be the random graph with $np_n\rightarrow\infty$  and $f(x,y)\geq \epsilon$ for some constant $\epsilon\in(0,1)$. If $s_2\asymp n^2$,
then  the  R\'{e}nyi index has the asymptotic expression as
\begin{eqnarray}\label{onerenyi}
\hskip 3.5cm
 \mathcal{R}_{1}=1-e^{-r_n+O_P\left(\frac{1}{n\sqrt{p_n}}\right)},\hskip 1cm r_n=\frac{1}{n}\sum_{i=1}^n\frac{\mu_i}{np_n\lambda_{0,1}}\log\left(\frac{\mu_i}{np_n\lambda_{0,1}}\right),
\end{eqnarray}
where the error rate $\frac{1}{n\sqrt{p_n}}$ cannot be improved.
\end{Theorem}

Based on Theorem \ref{theorem:2}, $\mathcal{R}_{1}$ mainly depends on $r_n$. 
For the Erd\"{o}s-R\'{e}nyi random graph, that is, $f(x,y)\equiv1$, it is obvious that $\lambda_{0,1}=1$, $\mu_i=(n-1)p_n$ and hence $\mathcal{R}_{1}=o_P(1)$. For $f(x,y)=e^{-\kappa x}e^{-\kappa y}$ with a positive constant $\kappa$,  $\mu_i=np_ne^{-\frac{\kappa i}{n}}\lambda_{1,0}(1+o(1))\asymp np_n$, then
$s_2\asymp n^2$.
The assumption of Theorem \ref{theorem:2} are satisfied. Straightforward calculation yields
\begin{equation}\label{renyi2}
r_n=g(\kappa)+o(1),
\end{equation}
where
\[g(\kappa)=-1+\frac{\kappa}{e^{\kappa}-1}-\log\left(\frac{e^{\kappa}-1}{\kappa e^{\kappa}}\right).\]
Note that
$\lim_{\kappa\rightarrow0^+}g(\kappa)=0$, $
\lim_{\kappa\rightarrow\infty}g(\kappa)=\infty$. 
Hence larger $\kappa$ produces more heterogeneous random graph. This is consistent with the case $\alpha\neq1$.

The assumption that $f(x,y)\geq \epsilon$ for a constant $\epsilon\in(0,1)$ in Theorem \ref{theorem:1} and Theorem \ref{theorem:2} can be relaxed and replaced by less restrictive  assumptions. However, the alternative assumptions are difficult to state and interpret and would lead to more complex proofs. For simplicity, we do not pursue this relaxation.  In addition, Theorem  \ref{theorem:1} and Theorem \ref{theorem:2} hold for sparse networks, since they allow $p_n=o(1)$ as long as $np_n\rightarrow\infty$.

\subsection{The R\'{e}nyi index of a power-law random graph}\label{mainpower}

Empirical studies have
shown that many real-world networks are highly heterogeneous, that is, the degrees of nodes follow a power-law distribution (\cite{CSN09,N03,VHHK19}). This motivates us to study whether the R\'{e}nyi index  of power-law random graph is actually close to one.

Given a positive constant $\tau$, let $W$ be a random variable following a distribution with power-law tail as
\begin{equation}\label{powerlaw}
\mathbb{P}(W>x)=x^{-\tau},\ \ x\geq 1.
\end{equation}
This distribution has heavy tail and the $k$-th moment of $W$ exists if and only if $k<\tau$. The distribution  (\ref{powerlaw}) is widely used to define power-law random graphs (\cite{BDM06,JLN10,JLS19}). Given independent and identically distributed random variables $\omega_1,\dots, \omega_n$ from distribution (\ref{powerlaw}), define a power-law random graph $\mathcal{G}(n,\tau)$ as
\[\mathbb{P}(A_{ij}=1|W)=p\frac{\tilde{\omega}_i\tilde{\omega}_j}{n},\]
where $W=(\omega_1,\dots, \omega_n)$, $\tilde{\omega}_i=\min\{\omega_i,\sqrt{n}\}$, $p\in(0,1)$ is a constant and $A_{ij}\ (1\leq i<j\leq n)$ are  independent conditional on $W$. 

The random graph $\mathcal{G}(n,\tau)$ was first defined in \cite{BM05,BM06} and the order of large  cliques was studied there. The cutoff $\sqrt{n}$ in $\tilde{\omega}_i$ guarantees the edge appearance probability is less than 1. This cutoff is common in algorithm analysis and random graph theory (\cite{BCH20,CGL16,YXL21}). We focus on the interesting regime $\tau\in(1,2)$ as in literature (\cite{BDM06,JLN10,JLS19}). 
Note that the edges $A_{ij} (1\leq i<j\leq n)$ are not independent and higher moments of $\tilde{\omega}_i$ are not bounded.  It is more challenging to derive the limit of the R\'{e}nyi index $\mathcal{R}_{\alpha}$ of $\mathcal{G}(n,\tau)$ for arbitrary $\alpha>0$.  In this paper, we only study $\mathcal{R}_{2}$.

\begin{Theorem}\label{theorem:3} Let $\mathcal{G}(n,\tau)$ be the power-law random graph with $\tau\in(1,2)$. Then
\begin{equation}\label{nus}
\mathcal{R}_{2}=1-O_P\left(\frac{1}{n^{1-\frac{\tau}{2}}}\right),
\end{equation}
 where the rate $\frac{1}{n^{1-\frac{\tau}{2}}}$ cannot be improved.
\end{Theorem}

According to Theorem \ref{theorem:3}, the R\'{e}nyi index $\mathcal{R}_2$ of $\mathcal{G}(n,\tau)$ converges to one in probability at rate $n^{\frac{\tau}{2}-1}$.  This indicates $\mathcal{G}(n,\tau)$ is extremely heterogeneous, consistent with empirical observations (\cite{CSN09,N03,VHHK19}).
Note that nodes  of $\mathcal{G}(n,\tau)$ have the same expected degree $p\left(\mathbb{E}[\omega_1]\right)^2$. In this sense, it seems $\mathcal{G}(n,\tau)$ is  homogeneous as the Erd\"{o}s-R\'{e}nyi random graph. However, the correlation between $A_{ij}$ and $A_{ik}$\ $(1\leq i\leq n,j\neq k)$ and the power-law tail property of $W$ jointly make the degrees extremely different so that $\mathcal{R}_2=1+o_P(1)$. Theorem \ref{theorem:3} provides an alternative justification that power-law random graph can be used as a generative model of extremely heterogeneous networks.

To conclude this section, we comment that Theorem \ref{theorem:1}, Theorem \ref{theorem:2} and Theorem \ref{theorem:3} jointly provide a theoretical justification that the R\'{e}nyi index is a reasonable measure of heterogeneity of networks. For homogeneous network, the R\'{e}nyi index is asymptotically zero. For extremely heterogeneous network, the R\'{e}nyi index is asymptotically one. For moderately heterogeneous network, the R\'{e}nyi index resides between zero and one.

\section{Simulation} \label{simu}

In this section, we conduct simulation study to evaluate finite-sample performance of the R\'{e}nyi index.

In this simulation, 20 graphs were generated from each random graph model described in Section \ref{main}, and the R\'{e}nyi index of each graph was calculated with $\alpha\in\{0.5,1,2,2.5,3,10\}$. Then the mean and standard deviation (sd) of the R\'{e}nyi indexes were computed, as well as the limit specified in Theorem \ref{theorem:1} or Theorem \ref{theorem:2}.

Firstly, we consider the heterogeneous Erd\"{o}s-R\'{e}nyi random graph $\mathcal{G}(n,p_n, f)$ with $f(x,y)=e^{-\kappa x}e^{-\kappa y}$ for a positive constant $\kappa$. The limit of the R\'{e}nyi index has a closed form given in (\ref{lrenyi}) for $\alpha\neq1$ and (\ref{renyi2}) for $\alpha=1$. With a little abuse of notation, we denote the limit as $\mathcal{R}_{\alpha}$. The model parameters we used to generate graphs, $\mathcal{R}_{\alpha}$, and the mean and standard deviation of the R\'{e}nyi indexes are listed in Table \ref{tab1}, Table \ref{tab2}, Table \ref{tab3}. 
As $n$ increases, the mean gets closer to the limit $\mathcal{R}_{\alpha}$, and $p_n$ highly affects the convergence speed. These findings coincide with the results in Theorem \ref{theorem:1} and Theorem \ref{theorem:2}. For homogeneous model ($\kappa=0.1$), the mean and limit $\mathcal{R}_{\alpha}$ almost vanish, while for heterogeneous model ($\kappa=25$) both are pretty large (greater than 0.8). This confirms that the R\'{e}nyi index can effectively measure heterogeneity of networks. In addition, the R\'{e}nyi indexes increase as $\alpha$ increases.

Now we consider the power-law random graph in Section \ref{mainpower}. The means and standard deviations (in parentheses) are summarized in Table \ref{tab4}. We point out that although the rate $\frac{1}{n^{1-\frac{\tau}{2}}}$ in Theorem \ref{theorem:3} only depends on $\tau$, the term $O_P\left(\frac{1}{n^{1-\frac{\tau}{2}}}\right)$ does involve constant $p$ in a complex way (see proof of Theorem \ref{theorem:3}). As a result, the values of  $p,\tau$ may significantly affect how close is the mean to the limit $\mathcal{R}_{2}=1$ in finite-sample case. Table \ref{tab4} shows all the means of the R\'{e}nyi indexes are larger than 0.6, indicating the power-law random graph is indeed heterogeneous. When $n=10,000$, most of the means are close to or larger than 0.90.

\begin{table}[h]
\begin{center}
\begin{tabular}{ |c||c|c| c|c|c|| c|c|c|} 
 \hline
 $n$ & $p_n$ &$\kappa$ & $\mathcal{R}_{0.5}$ & $mean$ & $sd$   & $\mathcal{R}_{1}$ & $mean$ & $sd$  \\  
  \hline
  \hline
 100  & 0.1 &0.1 & 0.0002&  0.0263&0.0039&0.0004 &0.0474& 0.0072 \\ 
 500  & 0.1 &0.1 & 0.0002& 0.0053& 0.0003 &0.0004&0.0103& 0.0005 \\ 
 2000  & 0.1 &0.1 &0.0002 & 0.0014 &0.0001 & 0.0004& 0.0029&0.0001    \\ 
 10000  & 0.1 &0.1 & 0.0002 &0.0004  & 0.0001   & 0.0004  & 0.0009 &0.0001      \\  
  \hline
 100  & 0.1 &4 & 0.238&0.6626 & 0.0326& 0.380& 0.7016 & 0.0280  \\ 
 500  & 0.1 &4 & 0.238& 0.3983& 0.0115&0.380 &0.4725 &  0.0098\\  
 2000  & 0.1 &4 &  0.238& 0.2780 &0.0035& 0.380&0.4069& 0.0028 \\  
  10000  & 0.1 &4  &0.238 & 0.2444 & 0.0007& 0.380 & 0.3861 & 0.0005 \\  
  \hline
 100  & 0.1 &25 &0.840&0.9763  &0.0073& 0.891 & 0.9913 &0.0109  \\ 
 500  & 0.1  &25 & 0.840& 0.9521&0.0072& 0.891&0.9533& 0.0080\\ 
 2000  & 0.1  &25 &0.840& 0.9078& 0.0025 &0.891&0.9198&0.0029 \\ 
  10000  & 0.1 & 25 & 0.840 & 0.8523 &0.0002 & 0.891 &0.8992 &  0.0004 \\  
 \hline
 \hline
 100  & 0.5 &0.1 &0.0002& 0.0032& 0.0005 &0.0004& 0.0063& 0.0008 \\ 
 500  & 0.5 &0.1 &0.0002&0.0008& 0.0001  &0.0004&0.0016& 0.0001\\  
 2000  & 0.5 &0.1 &0.0002&0.0003&0.0001   & 0.0004& 0.0007& 0.0001    \\  
 10000  & 0.5 &0.1 &0.0002& 0.0002 &0.0001& 0.0004& 0.0004  &  0.0001\\  
  \hline
 100  & 0.5 &4 &0.238&0.4029 &0.0413 & 0.380&0.4643& 0.0224 \\  
 500  & 0.5 &4 &0.238  & 0.2667 &0.0076& 0.380& 0.3984 &0.0065 \\ 
 2000  & 0.5 &4 & 0.238 &0.2445 &0.0008&0.380 &0.3854& 0.0011\\ 
 10000  & 0.5 &4 & 0.238 &0.2394 &0.0001&0.380  & 0.3817 & 0.0001 \\ 
  \hline
 100  & 0.5 &25 & 0.840 & 0.9583& 0.0185 &  0.891  &0.9573   & 0.0157 \\ 
 500  & 0.5  &25 &0.840&0.9024& 0.0072& 0.891   & 0.9045   &  0.0041\\ 
 2000  & 0.5  &25 & 0.840 & 0.8739& 0.0025 &  0.891  &  0.8998 & 0.0013\\ 
 10000  & 0.5 &25 &0.840 &  0.8443 & 0.0003& 0.891 & 0.8932 & 0.0002\\  
 \hline
\end{tabular}
\caption{The R\'{e}nyi index of heterogeneous Erd\"{o}s-R\'{e}nyi random graph with $\alpha=0.5, 1$. }\label{tab1}
\end{center}
\end{table}

\begin{table}[h]
\begin{center}
\begin{tabular}{ |c||c|c| c|c|c|| c|c|c|} 
 \hline
 $n$ & $p_n$ &$\kappa$ & $\mathcal{R}_{2}$ & $mean$ & $sd$   & $\mathcal{R}_{2.5}$ & $mean$ & $sd$  \\  
  \hline
  \hline
 100  & 0.1 &0.1 & 0.0008& 0.0916& 0.0113& 0.0010&   0.1082 & 0.0150\\ 
 500  & 0.1 &0.1 & 0.0008& 0.0208& 0.0014& 0.0010&  0.0251  & 0.0019 \\ 
 2000  & 0.1 &0.1 & 0.0008& 0.0058& 0.0001 &  0.0010 &  0.0072 & 0.0002 \\ 
 10000  & 0.1 &0.1 & 0.0008 &0.0018  & 0.0001   & 0.0010  &  0.0022 & 0.0001 \\
  \hline
 100  & 0.1 &4 & 0.517 & 0.7352&0.0412&0.553    & 0.7475    &  0.0386   \\ 
 500  & 0.1 &4 &0.517  &  0.5848& 0.0111& 0.553   & 0.6106 & 0.0099 \\  
 2000  & 0.1 &4 & 0.517&0.5361& 0.0029& 0.553 & 0.5714       & 0.0032    \\  
 10000  & 0.1 &4 & 0.517  & 0.5212  & 0.0006 &0.553 &0.5569&0.0006\\
  \hline
 100  & 0.1 &25 & 0.920 & 0.9761 & 0.0075& 0.926   & 0.9771 &   0.0075  \\ 
 500  & 0.1  &25 & 0.920 & 0.9572& 0.0064 &0.926   &0.9645   &  0.0054   \\ 
 2000  & 0.1  &25 &0.920  & 0.9344 &0.0032 & 0.926   &  0.9416   &  0.0034     \\ 
 10000  & 0.1 &25 & 0.920  & 0.9235  & 0.0004 & 0.926 &0.9300& 0.0003\\
 \hline
 \hline
 100  & 0.5 &0.1 & 0.0008&0.0126& 0.0020& 0.0010 &0.0158    & 0.0022\\ 
 500  & 0.5 &0.1 & 0.0008& 0.0032&0.0001 &0.0010&   0.0040& 0.0003\\  
 2000  & 0.5 &0.1 & 0.0008& 0.0014& 0.0001   &0.0010& 0.0017& 0.0001 \\  
 10000  & 0.5 &0.1 & 0.0008  &0.0009 & 0.0001     & 0.0010  & 0.0011  & 0.0001 \\
  \hline
 100  & 0.5 &4 & 0.517 &  0.5679 & 0.0253  & 0.553 & 0.6049& 0.0156 \\  
 500  & 0.5 &4 & 0.517 & 0.5315 &0.0051 & 0.553 & 0.5650   & 0.0057 \\ 
 2000  & 0.5 &4 & 0.517 &0.5211 & 0.0010& 0.553& 0.5563  & 0.0013 \\ 
 10000  & 0.5 &4 & 0.517  & 0.5187 &  0.0003 & 0.553  &0.5542 & 0.0002\\
  \hline
 100  & 0.5 &25 & 0.920 &  0.958  & 0.0170   &0.926 & 0.9612&0.0149\\ 
 500  & 0.5  &25 & 0.920  & 0.931   & 0.0052   &0.926  & 0.9347& 0.0055 \\ 
 2000  & 0.5  &25 & 0.920 & 0.923    &  0.0015  & 0.926   & 0.9290& 0.0015\\ 
 10000  & 0.5 &25 & 0.920& 0.921 &  0.0001 &0.926 &0.9270& 0.0002\\
 \hline
\end{tabular}
\caption{The R\'{e}nyi index of heterogeneous Erd\"{o}s-R\'{e}nyi random graph with $\alpha=2, 2.5$.}\label{tab2}
\end{center}
\end{table}

\begin{table}[h]
\begin{center}
\begin{tabular}{ |c||c|c| c|c|c|| c|c|c|} 
 \hline
 $n$ & $p_n$ &$\kappa$ & $\mathcal{R}_{3}$ & $mean$ & $sd$   & $\mathcal{R}_{10}$ & $mean$ & $sd$  \\  
  \hline
  \hline
 100  & 0.1 &0.1 & 0.0012& 0.1285 &0.0158 &0.0041 & 0.2612 & 0.0425   \\ 
 500  & 0.1 &0.1 & 0.0012& 0.0304& 0.0021   & 0.0041&0.0875 & 0.0059    \\ 
 2000  & 0.1 &0.1 &0.0012& 0.0086 & 0.0002 &0.0041 &0.0275 &0.0007    \\ 
   10000  & 0.1  & 0.1 & 0.0012 &0.0027 & 0.0001    & 0.0041  & 0.0091  & 0.0002 \\
  \hline
 100  & 0.1 &4 & 0.578 & 0.7540&0.0341 & 0.683& 0.8391 & 0.0312\\ 
 500  & 0.1 &4 &0.578  & 0.6400& 0.0090 & 0.683& 0.7603& 0.0171 \\  
 2000  & 0.1 &4 & 0.578 &0.5966 & 0.0042 & 0.683& 0.7112 & 0.0042 \\  
   10000  &0.1   & 4 &0.578 & 0.5820  & 0.0010    & 0.683  &0.6896 & 0.0012\\
  \hline
 100  & 0.1 &25 &0.930  &0.9686 & 0.0102 & 0.948& 0.9774 & 0.0065 \\ 
 500  & 0.1  &25 & 0.930 &0.9640 &0.0079 &0.948& 0.9720 &0.0064 \\ 
 2000  & 0.1  &25 & 0.930 & 0.9456&0.0026 & 0.948   &0.9643&  0.0027   \\ 
   10000  &  0.1 & 25  & 0.930  &0.9338 & 0.0005 & 0.948 &  0.9537 & 0.0010\\
 \hline
 \hline
 100  & 0.5 &0.1 & 0.0012&0.0203 & 0.0021 & 0.0041&0.0561 &  0.0076 \\ 
 500  & 0.5 &0.1 & 0.0012&0.0048 &0.0003& 0.0041&0.0154& 0.0008 \\  
 2000  & 0.5 &0.1 & 0.0012&0.0021&0.0001 &0.0041& 0.0070&0.0002 \\  
  10000  &   0.5 &0.1 & 0.0012  & 0.0014  &  0.0001    & 0.0041  & 0.0047 & 0.0001\\
  \hline
 100  & 0.5 &4 &0.578& 0.6178 &0.0202& 0.683 & 0.7265 & 0.0252 \\  
 500  & 0.5 &4 &0.578&0.5899 &0.0042& 0.683& 0.7028 & 0.0066\\ 
 2000  & 0.5 &4 &0.578& 0.5812 & 0.0012&0.683 & 0.6888 &0.0018\\ 
   10000  &  0.5 &4 &0.578  & 0.5795 &0.0002 & 0.683 &0.6847 & 0.0005\\
  \hline
 100  & 0.5 &25 & 0.930& 0.9618&0.0132 & 0.948 & 0.9661 &  0.0115   \\ 
 500  & 0.5  &25 & 0.930& 0.9394& 0.0050& 0.948 & 0.9585& 0.0047 \\ 
 2000  & 0.5  &25 & 0.930& 0.9337 & 0.0015& 0.948 & 0.9528 & 0.0017   \\ 
   10000  &  0.5  &25 & 0.930   & 0.9312 & 0.0001  & 0.948  & 0.9499 &0.0003 \\
 \hline
\end{tabular}
\caption{The R\'{e}nyi index of heterogeneous Erd\"{o}s-R\'{e}nyi random graph with $\alpha=3, 10$.}\label{tab3}
\end{center}
\end{table}

\begin{table}[h]
\begin{center}
\begin{tabular}{ |c|c|c| c|c|} 
 \hline
 $n$ & $p$ & $\tau=1.05$ & $\tau=1.50$ & $\tau=1.95$   \\  
  \hline
  \hline
        & 0.01 & 0.886(0.019)  & 0.939(0.013)  & 0.961(0.011) \\ 
 500  & 0.05 & 0.744(0.017)  & 0.811(0.021) & 0.855(0.017)\\ 
        & 0.25  & 0.645(0.013) & 0.656 (0.025) & 0.653(0.033)  \\ 
        & 0.95   & 0.622(0.016)  & 0.611(0.023)  & 0.541(0.044) \\
  \hline
  \hline
        & 0.01 & 0.877(0.007)  & 0.939(0.009)  & 0.961(0.008) \\ 
 1000  & 0.05 & 0.766(0.012) & 0.826(0.018)  & 0.855(0.018)  \\ 
        & 0.25  & 0.702(0.010)  & 0.691(0.024) & 0.671(0.026)  \\ 
        & 0.95  & 0.686(0.006) & 0.662(0.026)  & 0.549(0.064)\\
  \hline
  \hline
        & 0.01 & 0.886(0.009)  & 0.941(0.006) & 0.964(0.006) \\ 
 2000  & 0.05 & 0.794 (0.011)  & 0.831(0.013) & 0.858(0.013)  \\ 
        & 0.25  & 0.746(0.008)  & 0.731(0.016)  & 0.692(0.028)  \\ 
        & 0.95  & 0.741(0.011)  & 0.703(0.023)  & 0.606(0.032) \\
  \hline
\hline
        & 0.01 &0.927(0.002)  &0.946(0.002) & 0.964 (0.002) \\ 
 10000  & 0.05 &0.901(0.003)  &0.911(0.010) &  0.903 (0.010) \\ 
        & 0.25  &0.898(0.003)  &0.892(0.014) &0.851 (0.020) \\
        & 0.95  &0.895(0.001)  &0.894(0.011) & 0.825 (0.043) \\
  \hline
\end{tabular}
\caption{The R\'{e}nyi index of power-law random graph with $\alpha=2$.}\label{tab4}
\end{center}
\end{table}

\section{Proof of main results}\label{proof}
In this section, we provide detailed proof of the main results. Note that $\mathcal{R}_{\alpha}$ is a non-linear function of degrees $d_i$ as given in (\ref{renyiindex}). The degrees $d_i$ are not independent and may not be identically distributed. It is not feasible to directly apply the classical  Law of large number or Central limit theorem to get the limit of $\mathcal{R}_{\alpha}$. To overcome this issue, our strategy is to adopt the Taylor expansion to express the non-linear function of $d_i$ as a sum of polynomials of $d_i$ plus a remainder term. Then we carefully bound the remainder term and identify the limit and exact order of the polynomial terms. For convenience, let
\[\gamma_{k,l}=\frac{\sum_{i\neq j}f_i^{k}f_j^{k}f_{ij}^l}{n^2}.\]

\medskip

\noindent
{\bf Proof of Theorem \ref{theorem:1}:} The main challenge is that the degrees $d_i (1\leq i\leq n)$ are not independent and may not be identically distributed. The classical tools such as the law of large number and the central limit theorem can not be directly applied to derive the limit of $\mathcal{R}_{\alpha}$. Our proof strategy is: (a) use the Taylor expansion to expand $\sum_{i=1}^n\left(\frac{d_i}{n}\right)^{\alpha}$ at $\sum_{i=1}^n\left(\frac{\mu_i}{n}\right)^{\alpha}$ as a sum of polynomials in $d_i$ and a reminder term; (b) find the exact  order of the polynomial terms; (c) show the reminder term is bounded by the polynomial terms.  The key step is (c). We will use a truncation technique to control the reminder term.

To fix the idea, we consider the case $\alpha\in(0,3]\backslash \{1\}$ first.  Let $\mu_i=\mathbb{E}(d_i)$.
By Taylor expansion, we have
\begin{eqnarray}\nonumber
\left(\frac{d_i}{n}\right)^{\alpha}-\left(\frac{\mu_i}{n}\right)^{\alpha}&=&\alpha\left(\frac{\mu_i}{n}\right)^{\alpha-1}\left(\frac{d_i-\mu_i}{n}\right)+\frac{\alpha(\alpha-1)}{2!}\left(\frac{\mu_i}{n}\right)^{\alpha-2}\left(\frac{d_i-\mu_i}{n}\right)^2\\ \label{texpan}
&&+\frac{\alpha(\alpha-1)(\alpha-2)}{3!}X_{n,i}^{\alpha-3}\left(\frac{d_i-\mu_i}{n}\right)^3.
\end{eqnarray}
where $X_{n,i}$ is a random variable between $\frac{d_i}{n}$ and $\frac{\mu_i}{n}$. Summing both sides of (\ref{texpan}) over $i\in[n]$ yields
\begin{eqnarray}\nonumber
\sum_{i=1}^n\left(\frac{d_i}{n}\right)^{\alpha}-\sum_{i=1}^n\left(\frac{\mu_i}{n}\right)^{\alpha}&=&\alpha\sum_{i=1}^n\left(\frac{\mu_i}{n}\right)^{\alpha-1}\left(\frac{d_i-\mu_i}{n}\right)+\frac{\alpha(\alpha-1)}{2!}\sum_{i=1}^n\left(\frac{\mu_i}{n}\right)^{\alpha-2}\left(\frac{d_i-\mu_i}{n}\right)^2\\ \label{texpan00}
&&+\frac{\alpha(\alpha-1)(\alpha-2)}{3!}\sum_{i=1}^nX_{n,i}^{\alpha-3}\left(\frac{d_i-\mu_i}{n}\right)^3.
\end{eqnarray}
Next, we shall find the order of each term in the right-hand side of (\ref{texpan00}). We begin with the first term. For given $i\in[n]$, simple algebra yields
\begin{eqnarray*}
\mu_i=\sum_{j\neq i}\mathbb{E}(A_{ij})&=&\sum_{j\neq i}p_nf\left(\frac{i}{n},\frac{j}{n}\right)=np_nf_i,\\
\sum_{i=1}^n\left(\frac{\mu_i}{n}\right)^{\alpha}&=&p_n^{\alpha}\sum_{i=1}^nf_i^{\alpha}.
\end{eqnarray*}
Then
\begin{eqnarray}\nonumber
\sum_{i=1}^n\left(\frac{\mu_i}{n}\right)^{\alpha-1}\left(\frac{d_i-\mu_i}{n}\right)&=&p_n^{\alpha-1}\sum_{i=1}^nf_i^{\alpha-1}\frac{\sum_{j\neq i}(A_{ij}-f_{ij}p_n)}{n}\\ \label{new1}
&=&p_n^{\alpha-1}\frac{\sum_{i<j}(f_i^{\alpha-1}+f_j^{\alpha-1})(A_{ij}-f_{ij}p_n)}{n}.
\end{eqnarray}
Since $A_{ij}, (1\leq i<j\leq n)$ are independent and $\mathbb{E}[A_{ij}-f_{ij}p_n]=0$, then by (\ref{new1}) one has
\begin{eqnarray}\nonumber
&&\mathbb{E}\left[\sum_{i=1}^n\left(\frac{\mu_i}{n}\right)^{\alpha-1}\left(\frac{d_i-\mu_i}{n}\right)\right]^2\\
&=&p_n^{2\alpha-2}\mathbb{E}\left[\frac{\sum_{i<j}(f_i^{\alpha-1}+f_j^{\alpha-1})(A_{ij}-f_{ij}p_n)}{n}\right]^2\\ \nonumber
&=&p_n^{2\alpha-2}\frac{\sum_{i<j}\mathbb{E}(f_i^{\alpha-1}+f_j^{\alpha-1})^2(A_{ij}-f_{ij}p_n)^2}{n^2}\\ \nonumber
&=&p_n^{2\alpha-2}\left(\frac{\sum_{i\neq j}(f_i^{\alpha-1}+f_j^{\alpha-1})^2f_{ij}p_n}{2n^2}-\frac{\sum_{i\neq j}(f_i^{\alpha-1}+f_j^{\alpha-1})^2f_{ij}^2p_n^2}{2n^2}\right)\\ \nonumber
&=&p_n^{2\alpha-1}\frac{\sum_{i\neq j}f_i^{2\alpha-2}f_{ij}+\sum_{i\neq j}f_i^{\alpha-1}f_j^{\alpha-1}f_{ij}}{n^2}\\ \nonumber
&&-p_n^{2\alpha-1}\frac{p_n\sum_{i\neq j}f_i^{2\alpha-2}f_{ij}^2+p_n\sum_{i\neq j}f_i^{\alpha-1}f_j^{\alpha-1}f_{ij}^2}{n^2}\\ \label{inh1}
&=&p_n^{2\alpha-1}\left[\left(\lambda_{2\alpha-2,1}+\gamma_{\alpha-1,1}\right)-p_n\left(\lambda_{2\alpha-2,2}+\gamma_{\alpha-1,2}\right)\right]. 
\end{eqnarray}
 Since $f(x; y) \geq\epsilon > 0$, then $\lambda_{2\alpha-2,1}\asymp 1$, $\gamma_{\alpha-1,1}\asymp 1$, $\lambda_{2\alpha-2,2}\asymp 1$, $\gamma_{\alpha-1,2}\asymp 1$.
Hence the first term in the right-hand side of (\ref{texpan00}) is bounded by order $p_n^{\alpha-1}\sqrt{p_n}$. By Lemma \ref{lem1}, this is the exact order.

Secondly, we find the order of the second term in the right-hand side of (\ref{texpan00}). Note that
\begin{eqnarray}\nonumber
&&\sum_{i=1}^n\left(\frac{\mu_i}{n}\right)^{\alpha-2}\left(\frac{d_i-\mu_i}{n}\right)^2\\ \nonumber
&=&\frac{1}{n^2}\sum_{i=1}^np_n^{\alpha-2}f_i^{\alpha-2}\sum_{j,k\neq i}(A_{ij}-p_nf_{ij})(A_{ik}-p_nf_{ik})\\ \nonumber
&=&p_n^{\alpha-2}\frac{1}{n^2}\sum_{i\neq j}(A_{ij}-p_nf_{ij})^2f_i^{\alpha-2}+p_n^{\alpha-2}\frac{1}{n^2}\sum_{i\neq j\neq k}(A_{ij}-p_nf_{ij})(A_{ik}-p_nf_{ik})f_i^{\alpha-2}\\ \label{in3}
&=&S_1+S_2. 
\end{eqnarray}
We claim $S_2=o_P(S_1)$. To this end, we compute the second moment of $S_2$ and the first moment of $S_1$. Straightforward calculations yield
\begin{eqnarray}\nonumber
\mathbb{E}\left[S_1\right]&=&p_n^{\alpha-2}\frac{1}{n^2}\sum_{i\neq j}\mathbb{E}(A_{ij}-p_nf_{ij})^2f_i^{\alpha-2}\\ \nonumber
&=&p_n^{\alpha-2}\left(\frac{1}{n^2}\sum_{i\neq j}p_nf_{ij}f_i^{\alpha-2}-\frac{1}{n^2}\sum_{i\neq j}p_n^2f_{ij}^2f_i^{\alpha-2}\right)\\ \label{in2}
&=&p_n^{\alpha-1}\left(\lambda_{\alpha-2,1}-p_n\lambda_{\alpha-2,2}\right).
\end{eqnarray}
Note that $\lambda_{\alpha-2,1}\asymp 1$,  $\lambda_{\alpha-2,2}\asymp 1$, due to the assumption $f(x; y) \geq\epsilon > 0$. Then $S_1$ is bounded by order $p_n^{\alpha-1}$. By Lemma \ref{lem1}, this is the exact order.

Since $0\leq f_i\leq 1$ and $0\leq f_{ij}\leq 1$, then
\begin{eqnarray}\nonumber
\mathbb{E}\left[S_2^2\right]&\leq&\frac{p_n^{2(\alpha-2)}}{n^4}\sum_{\substack{i\neq j\neq k\\ r\neq s\neq t}}\mathbb{E}(A_{ij}-p_nf_{ij})(A_{ik}-p_nf_{ik})(A_{rs}-p_nf_{rs})(A_{rt}-p_nf_{rt})\\ \nonumber
&=&\frac{p_n^{2(\alpha-2)}}{n^4}O\left(\sum_{\substack{i\neq j\neq k}}\mathbb{E}(A_{ij}-p_nf_{ij})^2(A_{ik}-p_nf_{ik})^2\right)\\ \nonumber
&=&\frac{p_n^{2(\alpha-2)}}{n^4}O\left(\sum_{\substack{i\neq j\neq k}}p_n^2f_{ij}f_{ik}\right)\\ \label{in1}
&=&O\left(\frac{p_n^{2\alpha-2}}{n}\right).
\end{eqnarray}
Then $S_2=O_P\left(\frac{p_n^{\alpha-1}}{\sqrt{n}}\right)$. Hence the exact order of the second term in the right-hand side of (\ref{texpan00}) is $p_n^{\alpha-1}$ and $S_1$ is the leading term.

Next, we show the third term of (\ref{texpan00}) is bounded by $p_n^{\alpha-1}O\left(\frac{1}{\sqrt{np_n}}\right)$. If $\alpha=2$, then the third term in (\ref{texpan00}) vanishes. The desired result holds trivially. We only need to focus on the cases $\alpha\neq 1, 2$. Note that $X_{n,i}\geq0$. Then
\begin{eqnarray}\label{tecef4}
\mathbb{E}\left[\left|\sum_{i=1}^nX_{n,i}^{\alpha-3}\left(\frac{d_i-\mu_i}{n}\right)^3\right|\right]\leq \sum_{i=1}^n\mathbb{E}\left[X_{n,i}^{\alpha-3}\left|\frac{d_i-\mu_i}{n}\right|^3\right].
\end{eqnarray}
We shall show the right-hand side of (\ref{tecef4}) is bounded by $p_n^{\alpha-1}O\left(\frac{1}{\sqrt{np_n}}\right)$. 
Consider first the case $\alpha=3$. In this case, the expansion in Equation (\ref{texpan})
holds with $X_{n;i} = 1$, so that the analysis of Equation (\ref{tecef4}) is simpler.
Since $X_{n,i}=1$ for $\alpha=3$. By the Cauchy-Schwarz inequality, we have
\begin{eqnarray}\nonumber
\sum_{i=1}^n\mathbb{E}\left[\left|\frac{d_i-\mu_i}{n}\right|^3\right]
&\leq& \sum_{i=1}^n\sqrt{\mathbb{E}\left[\left(\frac{d_i-\mu_i}{n}\right)^6\right]}\\ \nonumber
&=&\sum_{i=1}^n\sqrt{\frac{\sum_{j_1,j_2,\dots,j_6\neq i}\mathbb{E}(A_{ij_1}-p_nf_{ij_1})(A_{ij_2}-p_nf_{ij_2})\dots(A_{ij_6}-p_nf_{ij_6})}{n^6}}\\ \nonumber
&=&\sum_{i=1}^n\sqrt{\frac{15\sum_{j_1\neq j_2\neq j_3\neq i}\mathbb{E}(A_{ij_1}-p_nf_{ij_1})^2(A_{ij_2}-p_nf_{ij_2})^2(A_{ij_3}-p_nf_{ij_3})^2}{n^6}}\\ \nonumber
&&+\sum_{i=1}^n\sqrt{\frac{15\sum_{j_1\neq j_2\neq i}\mathbb{E}(A_{ij_1}-p_nf_{ij_1})^4(A_{ij_2}-p_nf_{ij_2})^2}{n^6}}\\ \nonumber
&&+\sum_{i=1}^n\sqrt{\frac{20\sum_{j_1\neq j_2\neq i}\mathbb{E}(A_{ij_1}-p_nf_{ij_1})^3(A_{ij_2}-p_nf_{ij_2})^3}{n^6}}\\ \nonumber
&&+\sum_{i=1}^n\sqrt{\frac{\sum_{j_1\neq i}\mathbb{E}(A_{ij_1}-p_nf_{ij_1})^6}{n^6}}\\ \nonumber
&=&O\left(n\frac{\sqrt{n^3p_n^3}+\sqrt{n^2p_n^2}+\sqrt{np_n}}{\sqrt{n^6}}\right)\\ \label{inh2} 
&=&p_n^2O\left(\frac{1}{\sqrt{np_n}}+\frac{1}{np_n}+\frac{1}{np_n\sqrt{np_n}}\right)=p_n^2O\left(\frac{1}{\sqrt{np_n}}\right).
\end{eqnarray}
Hence, for $\alpha=3$, it follows that
\begin{equation}\label{tecef11}
\mathbb{E}\left[\left|\sum_{i=1}^nX_{n,i}^{\alpha-3}\left(\frac{d_i-\mu_i}{n}\right)^3\right|\right]=p_n^{\alpha-1}O\left(\frac{1}{\sqrt{np_n}}\right).
\end{equation}

Next we assume $\alpha\in(0,3)$ and $\alpha\neq 1,2$. Let $\delta\in(0,1)$ be an arbitrary small constant. Note that
\begin{eqnarray}\nonumber
&&\mathbb{E}\left[\left|\sum_{i=1}^nX_{n,i}^{\alpha-3}\left(\frac{d_i-\mu_i}{n}\right)^3\right|\right]=\mathbb{E}\left[\left|\sum_{i=1}^nX_{n,i}^{\alpha-3}\left(\frac{d_i-\mu_i}{n}\right)^3\left(I[X_{n,i}\geq\delta\frac{\mu_i}{n}]+I[X_{n,i}<\delta\frac{\mu_i}{n}]\right)\right|\right]\\ \label{mod1}
&\leq&\mathbb{E}\left[\sum_{i=1}^nX_{n,i}^{\alpha-3}\left|\frac{d_i-\mu_i}{n}\right|^3I[X_{n,i}\geq\delta\frac{\mu_i}{n}]\right]+\mathbb{E}\left[\sum_{i=1}^nX_{n,i}^{\alpha-3}\left|\frac{d_i-\mu_i}{n}\right|^3I\big[X_{n,i}<\delta\frac{\mu_i}{n}\big]\right]
\end{eqnarray}

Note
that, when $\alpha<3$,  then $\alpha-3<0$. If  $X_{n,i}\geq \delta\frac{\mu_i}{n}$, then $X_{n,i}^{\alpha-3}\leq \left(\delta\frac{\mu_i}{n}\right)^{\alpha-3}\leq \delta^{\alpha-3}p_n^{\alpha-3}f_i^{\alpha-3}$. So, it is possible to use the same approach as for the
case $\alpha=3$. By (\ref{tecef4}) and a similar calculation as in (\ref{inh2}) , we get
\begin{equation}\label{tecef21}
\mathbb{E}\left[\left|\sum_{i=1}^nX_{n,i}^{\alpha-3}\left(\frac{d_i-\mu_i}{n}\right)^3I[X_{n,i}\geq\delta\frac{\mu_i}{n}]\right|\right]\leq \delta^{\alpha-3}p_n^{\alpha-3}  \mathbb{E}\left[\sum_{i=1}^n\left|\frac{d_i-\mu_i}{n}\right|^3f_i^{\alpha-3}\right] =p_n^{\alpha-1}O\left(\frac{1}{\sqrt{np_n}}\right).
\end{equation}
The difficult case is $X_{n,i}< \delta\frac{\mu_i}{n}$.
Suppose $X_{n,i}< \delta\frac{\mu_i}{n}$.
Recall that $\frac{d_i}{n}\leq X_{n,i}\leq \frac{\mu_i}{n}$ or $\frac{\mu_i}{n}\leq X_{n,i}\leq \frac{d_i}{n}$. Then $X_{n,i}< \delta\frac{\mu_i}{n}$ implies $\frac{d_i}{n}\leq X_{n,i}\leq\delta \frac{\mu_i}{n}$. In this case, $\frac{d_i}{\mu_{i}}\leq \delta$.
 Dividing both sides of (\ref{texpan}) by $\left(\frac{\mu_i}{n}\right)^{\alpha}$ yields
\begin{eqnarray}\nonumber
\left(\frac{d_i}{\mu_{i}}\right)^{\alpha}-1=\alpha\left(\frac{d_i}{\mu_{i}}-1\right)+\frac{\alpha(\alpha-1)}{2}\left(\frac{d_i}{\mu_{i}}-1\right)^2+\frac{\alpha(\alpha-1)(\alpha-2)}{6}\frac{X_{n,1}^{\alpha-3}}{\left(\frac{\mu_i}{n}\right)^{\alpha-3}}\left(\frac{d_i}{\mu_{i}}-1\right)^3,
\end{eqnarray}
from which it follows that
\begin{eqnarray}\nonumber
&&\frac{\alpha(\alpha-1)(\alpha-2)}{6}\frac{X_{n,i}^{\alpha-3}}{\left(\frac{\mu_i}{n}\right)^{\alpha-3}}\left(\frac{d_i}{\mu_{i}}-1\right)^3\\ \label{eq1}
&=&-\frac{(\alpha-1)(\alpha-2)}{2}+\left(\frac{d_i}{\mu_{i}}\right)^{\alpha}+\alpha(\alpha-2)\frac{d_i}{\mu_{i}}-\frac{\alpha(\alpha-1)}{2}\left(\frac{d_i}{\mu_{i}}\right)^2.
\end{eqnarray}
Note that $\frac{d_i}{\mu_{i}}\geq 0$.
For a fixed $\alpha$, there exists a sufficiently small constant $\delta>0$ such that if $\frac{d_i}{\mu_{i}}\leq \delta$ then the right-hand side of (\ref{eq1}) is bounded away from zero and infinity. This implies that $X_{n,i}^{\alpha-3}\leq C\left(\frac{\mu_i}{n}\right)^{\alpha-3}$ for some constant $C>0$ and $C$ is independent of $i\in[n]$. Then similar to (\ref{tecef21}), we have
\begin{equation}\label{mod2}
\mathbb{E}\left[\left|\sum_{i=1}^nX_{n,i}^{\alpha-3}\left(\frac{d_i-\mu_i}{n}\right)^3I[X_{n,i}<\delta\frac{\mu_i}{n}]\right|\right] =p_n^{\alpha-1}O\left(\frac{1}{\sqrt{np_n}}\right).
\end{equation}
According to (\ref{mod1}), (\ref{tecef21}) and (\ref{mod2}),  (\ref{tecef11}) holds for $\alpha\in(0,3)$ and $\alpha\neq 1,2$.

By (\ref{tecef4}), (\ref{tecef11}) and (\ref{tecef21}), it follows that
the third term of (\ref{texpan00}) is bounded by $p_n^{\alpha-1}O_P\left(\frac{1}{\sqrt{np_n}}\right)$. Then the first two terms are the leading terms. By  (\ref{texpan00}), we get
\begin{eqnarray}\nonumber
\sum_{i=1}^n\left(\frac{d_i}{n}\right)^{\alpha}-p_n^{\alpha}\sum_{i=1}^nf_i^{\alpha}&=&\alpha\sum_{i=1}^n\left(p_nf_i\right)^{\alpha-1}\left(\frac{d_i-\mu_i}{n}\right)+\frac{\alpha(\alpha-1)}{2!}p_n^{\alpha-2}\frac{1}{n^2}\sum_{i\neq j}(A_{ij}-p_nf_{ij})^2f_i^{\alpha-2}\\ \nonumber
&&+O_P\left(\frac{p_n^{\alpha-1}}{\sqrt{np_n}}+\frac{p_n^{\alpha-1}}{\sqrt{n}}\right)\\ \label{pan0}
&=&O_P\left(p_n^{\alpha-1}\sqrt{p_n}\right)+O_P\left(p_n^{\alpha-1}\right)+O_P\left(\frac{p_n^{\alpha-1}}{\sqrt{np_n}}+\frac{p_n^{\alpha-1}}{\sqrt{n}}\right).
\end{eqnarray}
By Lemma \ref{lem1}, the rates $O_P\left(p_n^{\alpha-1}\sqrt{p_n}\right)$ and $O_P\left(p_n^{\alpha-1}\right)$ in (\ref{pan0}) are optimal. 
 Besides, $\frac{d}{n}=p_n\lambda_{0,1}+O_P\left(\frac{\sqrt{p_n}}{n}\right)$ and the rate $\frac{\sqrt{p_n}}{n}$ cannot be improved according to Lemma \ref{lem1}. Then (\ref{renyiasymp}) follows for $\alpha\in(0,3]\backslash \{1\}$.

 Now assume $\alpha\in(k-1,k]$ for any fixed integer $k\geq4$. By Taylor expansion, we have
 \begin{eqnarray}\nonumber
\sum_{i=1}^n\left(\frac{d_i}{n}\right)^{\alpha}-\sum_{i=1}^n\left(\frac{\mu_i}{n}\right)^{\alpha}&=&\alpha\sum_{i=1}^n\left(\frac{\mu_i}{n}\right)^{\alpha-1}\left(\frac{d_i-\mu_i}{n}\right)+\frac{\alpha(\alpha-1)}{2!}\sum_{i=1}^n\left(\frac{\mu_i}{n}\right)^{\alpha-2}\left(\frac{d_i-\mu_i}{n}\right)^2\\ \label{gtexpan00}
&&+\dots+\frac{\alpha(\alpha-1)\dots(\alpha-k+1)}{k!}\sum_{i=1}^nX_{n,i}^{\alpha-k}\left(\frac{d_i-\mu_i}{n}\right)^k,
\end{eqnarray}
where $X_{n,i}$ is between $\frac{d_i}{n}$ and $\frac{\mu_i}{n}$.
 To complete the proof, it suffices to show the first two terms of (\ref{gtexpan00}) are the leading terms.  More specifically, we show  only the first two terms ``matter" among the first $k- 1$ terms. Then we show the remainder
term is negligible using a truncation argument, analogous to the one used for
the case $\alpha<3$.

 For integer $t$ with $3\leq t\leq k$, we have
 \begin{eqnarray}\nonumber
\mathbb{E}\left[\left|\sum_{i=1}^n\left(\frac{\mu_i}{n}\right)^{\alpha-t}\left(\frac{d_i-\mu_i}{n}\right)^{t}\right|\right]&\leq& \sum_{i=1}^n\left(\frac{\mu_i}{n}\right)^{\alpha-t}\mathbb{E}\left[\left|\frac{d_i-\mu_i}{n}\right|^{t}\right]\\ \nonumber
 &\leq&\sum_{i=1}^n\left(\frac{\mu_i}{n}\right)^{\alpha-t}\sqrt{\mathbb{E}\left[\left(\frac{d_i-\mu_i}{n}\right)^{2t}\right]}\\ \label{geq1}
 &=&O\left(\frac{p_n^{\alpha-t}}{n^{t-1}}\sqrt{\mathbb{E}\left[\left(d_i-\mu_i\right)^{2t}\right]}\right).
 \end{eqnarray}
 Note that
 \begin{equation}\nonumber
 \mathbb{E}\left[\left(d_i-\mu_i\right)^{2t}\right]=\mathbb{E}\left[\sum_{j_1,j_2,\dots,j_{2t}}(A_{ij_1}-p_nf_{ij_1})\dots (A_{ij_{2t}}-p_nf_{ij_{2t}})\right].\end{equation}
Since $A_{ij}$ and $A_{il}$ are independent if $j\neq l$, then $\mathbb{E}[(A_{ij}-p_nf_{ij})(A_{il}-p_nf_{il})]=0$ if $j\neq l$. If there exists an index $j_{s}$ such that $j_{s}$ is not equal to $j_{r}$ for any $r=1,2,\dots,s-1,s+1,\dots, 2t$, then
 \[\mathbb{E}\left[\sum_{j_1,j_2,\dots,j_{2t}}(A_{ij_1}-p_nf_{ij_1})\dots (A_{ij_{2t}}-p_nf_{ij_{2t}})\right]=0.\]
 Hence, each index $j_{s}$ must equal another index $j_{r}$ with $r\neq s$. Then  
 \begin{equation}\nonumber
 \mathbb{E}\left[\left(d_i-\mu_i\right)^{2t}\right]=\sum_{s=1}^t\sum_{j_1,j_2,\dots,j_{s}:distinct}\mathbb{E}\left[(A_{ij_1}-p_nf_{ij_1})^{\lambda_1}\dots (A_{ij_{s}}-p_nf_{ij_{s}})^{\lambda_s}\right],
 \end{equation}
 where $\lambda_r\geq2$ are integers and $\lambda_1+\lambda_2+\dots+\lambda_s=2t$. It is easy to verify that for $\lambda_r\geq2$ ($r=1,2,\dots,s$),
 \[\mathbb{E}\left[(A_{ij_r}-p_nf_{ij_r})^{\lambda_r}\right]=(1-p_nf_{ij_r})^{\lambda_r}p_nf_{ij_r}+(-p_nf_{ij_r})^{\lambda_r}(1-p_nf_{ij_r})=O(p_n).\]
 Then
 $\mathbb{E}\left[\left(d_i-\mu_i\right)^{2t}\right]=O\left(\sum_{s=1}^tn^sp_n^s\right)=O(n^tp_n^t)$. By (\ref{geq1}), one has
  \begin{eqnarray}\label{geq2}
\mathbb{E}\left[\left|\sum_{i=1}^n\left(\frac{\mu_i}{n}\right)^{\alpha-t}\left(\frac{d_i-\mu_i}{n}\right)^{t}\right|\right]
 &=&O\left(\frac{p_n^{\alpha-1}}{\left(np_n\right)^{\frac{t}{2}-1}}\right),\hskip 1cm 3\leq t\leq k.
 \end{eqnarray}
 If $X_{n,i}\leq \delta\frac{\mu_i}{n}$ for a small constant $\delta>0$, by a similar argument as in (\ref{eq1}), one can get $X_{n,i}^{\alpha-k}\leq M\left(\frac{\mu_i}{n}\right)^{\alpha-k}$ for a large constant $M>0$. By  (\ref{geq2}), (\ref{pan0}) holds. Then the proof is complete.

\qed

\medskip

\begin{Lemma}\label{lem1}
Under the assumption of Theorem \ref{theorem:1}, the following results are true. Then
\begin{eqnarray}\label{lemeq1}
\frac{\sum_{i=1}^n\left(p_nf_i\right)^{\alpha-1}\left(\frac{d_i-\mu_i}{n}\right)}{p_n^{\alpha-1}\sqrt{p_n(\lambda_{2\alpha-2,1}+\gamma_{\alpha-1,1})-p_n^2(\lambda_{2\alpha-2,2}+\gamma_{\alpha-1,2})}}&\Rightarrow&\mathcal{N}(0,1),
\end{eqnarray}
\begin{eqnarray}\label{lemeq2}
\frac{1}{n^2}\sum_{i\neq j}(A_{ij}-p_nf_{ij})^2f_i^{\alpha-2}&=&p_n(\lambda_{\alpha-2,1}-p_n\lambda_{\alpha-2,2})+o_P(1),
\end{eqnarray}
and
\begin{equation}\label{dd}
\sum_{i<j}\frac{A_{ij}-p_nf_{ij}}{n\sqrt{p_n\lambda_{0,1}+p_n^2\lambda_{0,2}}}\Rightarrow\mathcal{N}(0,1).
\end{equation}
\end{Lemma}

\medskip

\noindent
{\bf Proof of Lemma \ref{lem1}:} (I).  By (\ref{new1}) and (\ref{inh1}), we get 
\[\frac{\sum_{i=1}^n\left(p_nf_i\right)^{\alpha-1}\left(\frac{d_i-\mu_i}{n}\right)}{p_n^{\alpha-1}\sqrt{p_n(\lambda_{2\alpha-2,1}+\gamma_{\alpha-1,1})-p_n^2(\lambda_{2\alpha-2,2}+\gamma_{\alpha-1,2})}}=\sum_{i<j}\frac{(f_i^{\alpha-1}+f_j^{\alpha-1})(A_{ij}-f_{ij}p_n)}{n\sqrt{p_n(\lambda_{2\alpha-2,1}+\gamma_{\alpha-1,1})-p_n^2(\lambda_{2\alpha-2,2}+\gamma_{\alpha-1,2})}}.\]
Note that $A_{ij}, (1\leq i<j\leq n)$ are independent and $0<\epsilon\leq f(x,y)\leq 1$. Then $\lambda_{2\alpha-2,2}\asymp 1$, $\lambda_{2\alpha-2,1}\asymp 1$, $\gamma_{\alpha-1,1}\asymp 1$, $\gamma_{\alpha-1,2}\asymp 1$ and
\begin{eqnarray*}
&&\sum_{i<j}\mathbb{E}\left[\frac{(f_i^{\alpha-1}+f_j^{\alpha-1})(A_{ij}-f_{ij}p_n)}{n\sqrt{p_n(\lambda_{2\alpha-2,1}+\gamma_{\alpha-1,1})-p_n^2(\lambda_{2\alpha-2,2}+\gamma_{\alpha-1,2})}}\right]^4\\
&=&O\left(\frac{\sum_{i<j}(f_i^{\alpha-1}+f_j^{\alpha-1})^4f_{ij}}{n^4p_n}\right)=o(1).
\end{eqnarray*}
By the Lyapunov central limit theorem, (\ref{lemeq1}) holds.

(II). Note that
\begin{eqnarray*}
&&\mathbb{E}\left[\frac{1}{n^2}\sum_{i<j}(f_i^{\alpha-2}+f_j^{\alpha-2})\big[(A_{ij}-p_nf_{ij})^2-p_nf_{ij}(1-p_nf_{ij})\big]\right]^2\\
&=&\frac{1}{n^4}\sum_{i<j}(f_i^{\alpha-2}+f_j^{\alpha-2})^2\mathbb{E}\big[(A_{ij}-p_nf_{ij})^2-p_nf_{ij}(1-p_nf_{ij})\big]^2\\
&=&O\left(\frac{1}{n^4}\sum_{i<j}(f_i^{\alpha-2}+f_j^{\alpha-2})^2p_nf_{ij}\right)\\
&=&O\left(\frac{p_n}{n^4}\sum_{i\neq j}f_{ij}f_i^{2(\alpha-2)}+\frac{p_n}{n^4}\sum_{i\neq j}f_{ij}f_i^{\alpha-2}f_j^{\alpha-2}\right)=o(1).
\end{eqnarray*}
Hence (\ref{lemeq2}) holds.\\

(III). Note that
\begin{eqnarray*}
\mathbb{E}\left[\frac{\sum_{i<j}(A_{ij}-p_nf_{ij})}{n^2}\right]^2=\frac{\sum_{i<j}\mathbb{E}(A_{ij}-p_nf_{ij})^2}{n^4}=\frac{\sum_{i<j}(p_nf_{ij}-p_n^2f_{ij}^2)}{n^4}=\frac{p_n\lambda_{0,1}+p_n^2\lambda_{0,2}}{n^2}.
\end{eqnarray*}
Since $f(x,y)\geq\epsilon>0$, then $\lambda_{0,1}\asymp 1 $, $\lambda_{0,2}\asymp 1 $ and
\[\sum_{i<j}\frac{\mathbb{E}(A_{ij}-p_nf_{ij})^4}{\left(n\sqrt{p_n\lambda_{0,1}+p_n^2\lambda_{0,2}}\right)^4}=O\left(\frac{\sum_{i<j}f_{ij}}{n^4p_n}\right)=o(1).\]
By the Lyapunov central limit theorem, (\ref{dd}) holds.

\qed

\medskip

\noindent
{\bf Proof of Theorem \ref{theorem:2}:} The proof strategy is the same as the proof of Theorem \ref{theorem:1}.
Note that $\frac{d}{n}=p_n\lambda_{0,1}+O_P\left(\frac{\sqrt{p_n}}{n}\right)$ by Lemma \ref{lem1} and $d=\frac{1}{n}\sum_{i=1}^nd_i$. Then we have
\begin{eqnarray}\nonumber
\frac{1}{n}\sum_{i=1}^n\frac{d_i}{d}\log\frac{d_i}{d}&=&
\frac{1}{n}\sum_{i=1}^n\frac{d_i}{d}\log\left(\frac{np_n\lambda_{0,1}}{d}\frac{d_i}{np_n\lambda_{0,1}}\right)\\ \nonumber
&=&
\frac{1}{n}\sum_{i=1}^n\frac{d_i}{d}\log\left(\frac{np_n\lambda_{0,1}}{d}\right)+\frac{1}{n}\sum_{i=1}^n\frac{d_i}{d}\log\left(\frac{d_i}{np_n\lambda_{0,1}}\right)\\ \nonumber
&=&\log\left(\frac{np_n\lambda_{0,1}}{d}\right)+\frac{np_n\lambda_{0,1}}{d}\frac{1}{n}\sum_{i=1}^n\frac{d_i}{np_n\lambda_{0,1}}\log\frac{d_i}{np_n\lambda_{0,1}}\\ \label{did}
&=&O_P\left(\frac{1}{n\sqrt{p_n}}\right)+\frac{np_n\lambda_{0,1}}{d}\frac{1}{n}\sum_{i=1}^n\frac{d_i}{np_n\lambda_{0,1}}\log\frac{d_i}{np_n\lambda_{0,1}}.
\end{eqnarray}
It suffices to get the limit of $\sum_{i=1}^n\frac{d_i}{np_n\lambda_{0,1}}\log\frac{d_i}{np_n\lambda_{0,1}}$. Recall that $\mu_i=\mathbb{E}(d_i)=\sum_{j\neq i}p_nf_{ij}$. By the Taylor expansion, we have
\begin{eqnarray}\nonumber
\frac{d_i}{np_n\lambda_{0,1}}\log\left(\frac{d_i}{np_n\lambda_{0,1}}\right)&=&\frac{d_i}{np_n\lambda_{0,1}}\log\left(\frac{\mu_i}{np_n\lambda_{0,1}}\right)+\frac{d_i}{\mu_i}\left(\frac{d_i-\mu_i}{np_n\lambda_{0,1}}\right)\\ \label{aeq1}
&&-\frac{np_n\lambda_{0,1}d_i}{2\mu_i^2}\left(\frac{d_i-\mu_i}{np_n\lambda_{0,1}}\right)^2+\frac{1}{3X_{n,i}^3}\frac{d_i}{np_n\lambda_{0,1}}\left(\frac{d_i-\mu_i}{np_n\lambda_{0,1}}\right)^3,
\end{eqnarray}
where $\frac{d_i}{np_n\lambda_{0,1}}\leq X_{n,i}\leq \frac{\mu_i}{np_n\lambda_{0,1}}$ or $\frac{\mu_i}{np_n\lambda_{0,1}}\leq X_{n,i}\leq \frac{d_i}{np_n\lambda_{0,1}}$. Summing both sides of (\ref{aeq1}) over $i\in[n]$ yields
\begin{eqnarray}\nonumber
\sum_{i=1}^n\frac{d_i}{np_n\lambda_{0,1}}\log\left(\frac{d_i}{np_n\lambda_{0,1}}\right)&=&\sum_{i=1}^n\frac{d_i}{np_n\lambda_{0,1}}\log\left(\frac{\mu_i}{np_n\lambda_{0,1}}\right)+\sum_{i=1}^n\frac{d_i}{\mu_i}\left(\frac{d_i-\mu_i}{np_n\lambda_{0,1}}\right)\\  \nonumber
&&-\sum_{i=1}^n\frac{np_n\lambda_{0,1}d_i}{2\mu_i^2}\left(\frac{d_i-\mu_i}{np_n\lambda_{0,1}}\right)^2+\sum_{i=1}^n\frac{1}{3X_{n,i}^3}\frac{d_i}{np_n\lambda_{0,1}}\left(\frac{d_i-\mu_i}{np_n\lambda_{0,1}}\right)^3.\\ \label{aeq2}
\end{eqnarray}
Next we isolate the leading terms in the right-hand side of (\ref{aeq2}). More specifically, we show the first term is the leading term, and the second term, the third terms and the remainder term are of smaller order. For the remainder term, a truncation technique as in the proof of Theorem \ref{theorem:1} will be used.

 Firstly, we consider the second of (\ref{aeq2}). Note that
\begin{eqnarray}\label{ceq3}
\sum_{i=1}^n\frac{d_i}{\mu_i}\left(\frac{d_i-\mu_i}{np_n\lambda_{0,1}}\right)
=\sum_{i=1}^n\frac{(d_i-\mu_i)^2}{\mu_inp_n\lambda_{0,1}}+\sum_{i=1}^n\frac{d_i-\mu_i}{np_n\lambda_{0,1}}.
\end{eqnarray}
We find the order of each  term in the right-hand side of (\ref{ceq3}). Recall that $A_{ij},(1\leq i<j\leq n)$ are independent. Then straightforward calculations yield
\begin{eqnarray}\nonumber
\mathbb{E}\left[\sum_{i=1}^n\frac{(d_i-\mu_i)^2}{\mu_inp_n\lambda_{0,1}}\right]&=&\sum_{i=1}^n\frac{\sum_{j\neq i}\mathbb{E}(A_{ij}-p_nf_{ij})^2}{\mu_inp_n\lambda_{0,1}}\\ \nonumber
&=&\sum_{i=1}^n\frac{\sum_{j\neq i}p_nf_{ij}(1-p_nf_{ij})}{\mu_inp_n\lambda_{0,1}}\\ \nonumber
&=&\sum_{i=1}^n\frac{\sum_{j\neq i}p_nf_{ij}-\sum_{j\neq i}p_n^2f_{ij}^2}{\mu_inp_n\lambda_{0,1}}\\ \label{ceq1}
&=&\frac{1}{p_n\lambda_{0,1}}\left(1-p_n\frac{1}{n}\sum_{i=1}^n\frac{\sum_{j\neq i}f_{ij}^2}{\sum_{j\neq i}f_{ij}}\right),
\end{eqnarray}
and
\begin{eqnarray}\label{ceq2}
\mathbb{E}\left[\sum_{i=1}^n\frac{d_i-\mu_i}{np_n\lambda_{0,1}}\right]^2=\sum_{i\neq j}^n\frac{\mathbb{E}(A_{ij}-p_nf_{ij})^2}{n^2p_n^2\lambda_{0,1}^2}=\frac{\sum_{i\neq j}^np_nf_{ij}-\sum_{i\neq j}^np_n^2f_{ij}^2}{n^2p_n^2\lambda_{0,1}^2}=\frac{1}{p_n\lambda_{0,1}}-\frac{\lambda_{0,2}}{\lambda_{0,1}^2}.
\end{eqnarray}
Note that $\frac{1}{n}\sum_{i=1}^n\frac{\sum_{j\neq i}f_{ij}^2}{\sum_{j\neq i}f_{ij}}\leq 1$. Then (\ref{ceq3}) has order $O_P\left(\frac{1}{p_n}\right)$.

Next, we get the order of the third term in the right-hand side of (\ref{aeq2}).  Simple algebra yields
\begin{eqnarray}\label{ceq10}
\sum_{i=1}^n\frac{np_n\lambda_{0,1}d_i}{\mu_i^2}\left(\frac{d_i-\mu_i}{np_n\lambda_{0,1}}\right)^2=\frac{1}{np_n\lambda_{0,1}}\sum_{i=1}^n\frac{(d_i-\mu_i)^3}{\mu_i^2}+\frac{1}{np_n\lambda_{0,1}}\sum_{i=1}^n\frac{(d_i-\mu_i)^2}{\mu_i}.
\end{eqnarray}
Now we get an upper bound of the two terms in (\ref{ceq10}).
By the Cauchy-Schwarz inequality, one gets
\begin{eqnarray}\nonumber
&&\frac{1}{np_n\lambda_{0,1}}\mathbb{E}\left[\left|\sum_{i=1}^n\frac{(d_i-\mu_i)^3}{\mu_i^2}\right|\right]\\ \nonumber
&\leq&\frac{1}{np_n\lambda_{0,1}}\sum_{i=1}^n\mathbb{E}\left[\frac{|d_i-\mu_i|^3}{\mu_i^2}\right]\\ \nonumber
&\leq&\frac{1}{np_n\lambda_{0,1}}\sum_{i=1}^n\frac{1}{\mu_i^2}\sqrt{\mathbb{E}(d_i-\mu_i)^6}\\ \nonumber
&\leq&\frac{1}{np_n\lambda_{0,1}}\sum_{i=1}^n\frac{1}{\mu_i^2}\sqrt{\sum_{j_1\neq j_2\neq j_3\neq i}\mathbb{E}(A_{ij_1}-p_nf_{ij_1})^2(A_{ij_2}-p_nf_{ij_2})^2(A_{ij_3}-p_nf_{ij_3})^2}\\ \nonumber
&&+\frac{1}{np_n\lambda_{0,1}}\sum_{i=1}^n\frac{1}{\mu_i^2}\sqrt{\sum_{j_1\neq j_2\neq i}\mathbb{E}(A_{ij_1}-p_nf_{ij_1})^3(A_{ij_2}-p_nf_{ij_2})^3}\\ \nonumber
&&+\frac{1}{np_n\lambda_{0,1}}\sum_{i=1}^n\frac{1}{\mu_i^2}\sqrt{\sum_{j_1\neq j_2\neq i}\mathbb{E}(A_{ij_1}-p_nf_{ij_1})^4(A_{ij_2}-p_nf_{ij_2})^2}
\\ \nonumber
&&+\frac{1}{np_n\lambda_{0,1}}\sum_{i=1}^n\frac{1}{\mu_i^2}\sqrt{\sum_{j_1\neq i}\mathbb{E}(A_{ij_1}-p_nf_{ij_1})^6}\\ \nonumber
&=&\frac{1}{np_n\lambda_{0,1}}O\left(\sum_{i=1}^n\frac{1}{\sqrt{\mu_i}}+\sum_{i=1}^n\frac{1}{\mu_i}+\sum_{i=1}^n\frac{1}{\sqrt{\mu_i^3}}\right)\\ \label{ceq5}
&=&\frac{1}{p_n}O\left(\frac{1}{\sqrt{np_n}}\right).
\end{eqnarray}
 Then the first term of (\ref{ceq10}) is bounded by $\frac{1}{p_n}O_P\left(\frac{1}{\sqrt{np_n}}\right)$.  By (\ref{ceq1}), the second term is bounded by  $O_P\left(\frac{1}{p_n}\right)$.

Next, we consider the last term in the right-hand side of (\ref{aeq2}). Let $\delta\in(0,1)$ be an {\bf arbitrary} small constant.  We shall find an upper bound of the last term of (\ref{aeq2}) in two cases: $X_{n,i}\geq\delta\frac{\mu_i}{np_n\lambda_{0,1}}$ and $X_{n,i}< \delta\frac{\mu_i}{np_n\lambda_{0,1}}$.
If $X_{n,i}\geq \delta\frac{\mu_i}{np_n\lambda_{0,1}}$, then
\begin{equation}\label{ceq6}
\frac{1}{3X_{n,i}^3}\frac{d_i}{np_n\lambda_{0,1}}\left|\frac{d_i-\mu_i}{np_n\lambda_{0,1}}\right|^3\leq \frac{1}{3\delta^3}\frac{d_i}{np_n\lambda_{0,1}}\left|\frac{d_i-\mu_i}{\mu_i}\right|^3.
\end{equation}
Suppose $X_{n,i}< \delta\frac{\mu_i}{np_n\lambda_{0,1}}$. If $X_{n,i}<\frac{d_i}{np_n\lambda_{0,1}}$, then $X_{n,i}$ cannot be between $\frac{d_i}{np_n\lambda_{0,1}}$ and $\frac{\mu_i}{np_n\lambda_{0,1}}$.  Therefore, $\frac{d_i}{np_n\lambda_{0,1}}\leq X_{n,i}< \delta\frac{\mu_i}{np_n\lambda_{0,1}}$. Then $\frac{d_i}{\mu_i}\leq \delta$. Since $-\log x\rightarrow\infty$ as $x\rightarrow0^+$ and $\frac{d_i}{\mu_i}\geq0$, for small enough $\delta$, by (\ref{aeq1}) we have
\[\frac{\left(\frac{\mu_i}{np_n\lambda_{0,1}}\right)^3}{3X_{n,i}^3}=\frac{-\log\left(\frac{d_i}{\mu_i}\right)+\left(\frac{d_i}{\mu_i}-1\right)-\frac{1}{2}\left(\frac{d_i}{\mu_i}-1\right)^2}{(1-\frac{d_i}{\mu_i})^3}\leq -2\log\left(\frac{d_i}{\mu_i}\right).\]
Consequently, it follows that
\begin{equation*}
\frac{1}{3X_{n,i}^3}\frac{d_i}{np_n\lambda_{0,1}}\left|\frac{d_i-\mu_i}{np_n\lambda_{0,1}}\right|^3\leq -2\log\left(\frac{d_i}{\mu_i}\right)\frac{d_i}{\mu_i}\frac{|d_i-\mu_i|^3}{\mu_i^2np_n\lambda_{0,1}}.
\end{equation*}
Note that $\lim_{x\rightarrow0^+}x\log x=o(1)$. For small enough $\delta$, it follows that $-2\log\left(\frac{d_i}{\mu_i}\right)\frac{d_i}{\mu_i}\leq 1$ and hence
\begin{equation}\label{ceq7}
\frac{1}{3X_{n,i}^3}\frac{d_i}{np_n\lambda_{0,1}}\left|\frac{d_i-\mu_i}{np_n\lambda_{0,1}}\right|^3\leq \frac{|d_i-\mu_i|^3}{\mu_i^2np_n\lambda_{0,1}}.
\end{equation}

By (\ref{ceq6}) and  (\ref{ceq7}), for a fixed small constant $\delta\in(0,1)$, one has
\begin{eqnarray}\nonumber
\sum_{i=1}^n\frac{1}{3X_{n,i}^3}\frac{d_i}{np_n\lambda_{0,1}}\left|\frac{d_i-\mu_i}{np_n\lambda_{0,1}}\right|^3&=&\sum_{i=1}^n\frac{1}{3X_{n,i}^3}\frac{d_i}{np_n\lambda_{0,1}}\left|\frac{d_i-\mu_i}{np_n\lambda_{0,1}}\right|^3I[X_{n,i}<\delta\frac{\mu_i}{np_n\lambda_{0,1}}]\\ \nonumber
&&+ \sum_{i=1}^n\frac{1}{3X_{n,i}^3}\frac{d_i}{np_n\lambda_{0,1}}\left|\frac{d_i-\mu_i}{np_n\lambda_{0,1}}\right|^3I[X_{n,i}\geq\delta\frac{\mu_i}{np_n\lambda_{0,1}}]\\ \label{ceq9}
&\leq& \frac{1}{np_n\lambda_{0,1}}\sum_{i=1}^n\frac{|d_i-\mu_i|^3}{\mu_i^2}+\frac{1}{3\delta^3}\sum_{i=1}^n\frac{d_i}{np_n\lambda_{0,1}}\left|\frac{d_i-\mu_i}{\mu_i}\right|^3.
\end{eqnarray}
By  (\ref{ceq5}), it follows that
\begin{eqnarray}\nonumber
\frac{1}{np_n\lambda_{0,1}}\mathbb{E}\left[\left|\sum_{i=1}^n\frac{d_i}{\mu_i}\frac{(d_i-\mu_i)^3}{\mu_i^2}\right|\right]&\leq &\frac{1}{np_n\lambda_{0,1}}\sum_{i=1}^n\sqrt{\mathbb{E}\left(\frac{d_i}{\mu_i}\right)^2}\sqrt{\frac{\mathbb{E}(d_i-\mu_i)^6}{\mu_i^4}}\\ \nonumber
&\leq&\frac{1}{np_n\lambda_{0,1}}\sum_{i=1}^n(1+\frac{1}{\sqrt{\mu_i}})\sqrt{\frac{\mathbb{E}(d_i-\mu_i)^6}{\mu_i^4}}\\ \label{ceq8}
&=&\frac{1}{p_n}O\left(\frac{1}{\sqrt{np_n}}\right).
\end{eqnarray}

Hence by(\ref{aeq2}), (\ref{ceq3}), (\ref{ceq1}), (\ref{ceq2}),  (\ref{ceq5}), (\ref{ceq9}), (\ref{ceq8}) and (\ref{ceq10}), we get that
\begin{eqnarray}\label{che1}
&&\sum_{i=1}^n\frac{d_i}{np_n\lambda_{0,1}}\log\left(\frac{d_i}{np_n\lambda_{0,1}}\right)-\sum_{i=1}^n\frac{\mu_i}{np_n\lambda_{0,1}}\log\left(\frac{\mu_i}{np_n\lambda_{0,1}}\right)\\ \label{aeq12} \nonumber
&=&\sum_{i=1}^n\frac{d_i-\mu_i}{np_n\lambda_{0,1}}\left(1+\log\left(\frac{\mu_i}{np_n\lambda_{0,1}}\right)\right)
+\frac{1}{2np_n\lambda_{0,1}}\sum_{i=1}^n\frac{(d_i-\mu_i)^2}{\mu_i}+\frac{1}{p_n}O_P\left(\frac{1}{\sqrt{np_n}}\right). 
\end{eqnarray}
Further, it follows from Lemma \ref{lem2} that
\begin{eqnarray}\nonumber
&&\frac{1}{\sqrt{s_2p_n}}\sum_{i=1}^nd_i\log\left(\frac{d_i}{np_n\lambda_{0,1}}\right)-\frac{1}{\sqrt{s_2p_n}}\sum_{i=1}^n\mu_i\log\left(\frac{\mu_i}{np_n\lambda_{0,1}}\right)\\ \nonumber
&=&\sum_{i=1}^n\frac{d_i-\mu_i}{\sqrt{s_2p_n}}\left(1+l_i\right)+\frac{1}{2}\sum_{i=1}^n\frac{\sum_{j\neq i}(A_{ij}-p_nf_{ij})^2}{\mu_i\sqrt{s_2p_n}}+O_P\left(\frac{\sqrt{np_n}}{p_n\sqrt{s_2p_n}}+\frac{\sqrt{n}}{\sqrt{s_2p_n}}\right)\\ \label{oner1}
&=&O_P\left(1\right)+O_P\left(\frac{p_n\tau_1}{\sqrt{s_2p_n}}\right),
\end{eqnarray}
and the rates $O_P\left(1\right)$ and $O_P\left(\frac{p_n\tau_1}{\sqrt{p_ns_2}}\right)$ cannot be improved. Then (\ref{onerenyi}) follows from (\ref{did}) and (\ref{oner1}).

\qed

\begin{Lemma}\label{lem2}
Under the assumptions of Theorem \ref{theorem:2},
the following results are true.
\begin{eqnarray}\label{lemaeq1}
\sum_{i=1}^n\frac{d_i-\mu_i}{\sqrt{s_2p_n}}\left(1+l_i\right)&\Rightarrow&\mathcal{N}(0,1),\\ \label{lemaeq2}
\sum_{i=1}^n\frac{(d_i-\mu_i)^2}{\mu_i\sqrt{s_2p_n}}&=&\frac{p_n\tau_1-p_n^2\tau_{1,2}}{\sqrt{s_2p_n}}+o_P(1).
\end{eqnarray}
\end{Lemma}

\noindent
{\bf Proof of Lemma \ref{lem2}}. We firstly prove (\ref{lemaeq1}). Note that
\begin{eqnarray*}
\sum_{i=1}^n(d_i-\mu_i)(1+l_i)&=&\sum_{i<j}(A_{ij}-p_nf_{ij})(2+l_i+l_j).
\end{eqnarray*}
Then
\begin{eqnarray*}
\mathbb{E}\left[\sum_{i<j}(A_{ij}-p_nf_{ij})(2+l_i+l_j)\right]^2&=&\sum_{i<j}\mathbb{E}(A_{ij}-p_nf_{ij})^2(2+l_i+l_j)^2\\
&=&\sum_{i<j}(2+l_i+l_j)^2p_nf_{ij}(1-p_nf_{ij})=s_2p_n.
\end{eqnarray*}
Besides, 
\begin{eqnarray*}
\sum_{i<j}\frac{\mathbb{E}(A_{ij}-p_nf_{ij})^4(2+l_i+l_j)^4}{s_2^2p_n^2}=O\left(\frac{\sum_{i<j}(2+l_i+l_j)^4p_nf_{ij}}{s_2^2p_n^2}\right)=O\left(\frac{s_4}{s_2^2p_n}\right)=o(1).
\end{eqnarray*}
By the Lyapunov central limit theorem, we have
\[\frac{\sum_{i<j}(A_{ij}-p_nf_{ij})(2+l_i+l_j)}{\sqrt{s_2p_n}}\Rightarrow\mathcal{N}(0,1).\]

Next we  prove (\ref{lemaeq2}) .
Note that
\begin{eqnarray}\nonumber
\sum_{i=1}^n\frac{(d_i-\mu_i)^2}{\mu_i}&=&\sum_{i=1}^n\frac{\sum_{j\neq k\neq i}(A_{ij}-p_nf_{ij})(A_{ik}-p_nf_{ik})}{\mu_i}\\ \label{neq1}
&&+\sum_{i=1}^n\frac{\sum_{j\neq i}(A_{ij}-p_nf_{ij})^2}{\mu_i}.
\end{eqnarray}
Since
\begin{eqnarray*}
&&\mathbb{E}\left[\sum_{i=1}^n\frac{\sum_{j\neq k\neq i}(A_{ij}-p_nf_{ij})(A_{ik}-p_nf_{ik})}{\mu_i}\right]^2\\
&=&\sum_{i=1}^n\frac{\sum_{j\neq k\neq i}\mathbb{E}(A_{ij}-p_nf_{ij})^2(A_{ik}-p_nf_{ik})^2}{\mu_i^2}=O\left(n\right),
\end{eqnarray*}
then
\begin{eqnarray}\nonumber
\sum_{i=1}^n\frac{(d_i-\mu_i)^2}{\mu_i}&=&\sum_{i=1}^n\frac{\sum_{j\neq i}(A_{ij}-p_nf_{ij})^2}{\mu_i}+O_P\left(\sqrt{n}\right).
\end{eqnarray}
Note that
\[\sum_{i=1}^n\frac{\sum_{j\neq i}(A_{ij}-p_nf_{ij})^2}{\mu_i}=\sum_{i<j}\left(\frac{1}{\mu_i}+\frac{1}{\mu_j}\right)(A_{ij}-p_nf_{ij})^2,\]
\[\sum_{i<j}\left(\frac{1}{\mu_i}+\frac{1}{\mu_j}\right)\mathbb{E}(A_{ij}-p_nf_{ij})^2=\sum_{i<j}\left(\frac{1}{\mu_i}+\frac{1}{\mu_j}\right)p_nf_{ij}(1-p_nf_{ij}),\]
and
\[Var\left(\sum_{i<j}\left(\frac{1}{\mu_i}+\frac{1}{\mu_j}\right)(A_{ij}-p_nf_{ij})^2\right)=p_n\tau_2(1+o(1)).\]
Since $\tau_2=o(s_2)$, then

\[\frac{\sum_{i<j}\left(\frac{1}{\mu_i}+\frac{1}{\mu_j}\right)(A_{ij}-p_nf_{ij})^2}{\sqrt{s_2p_n}}=\frac{\sum_{i<j}\left(\frac{1}{\mu_i}+\frac{1}{\mu_j}\right)p_nf_{ij}(1-p_nf_{ij})}{\sqrt{s_2p_n}}+o_P(1).\]
\qed

\begin{Lemma}\label{lemma3}
For positive $k$ with $k\neq \tau$, we have 
\[\mathbb{E}(\tilde{\omega}_1^k)=n^{\frac{k-\tau}{2}}\frac{k}{k-\tau}-\frac{\tau}{k-\tau}.\]
\end{Lemma}
\noindent
\textit{Proof of Lemma \ref{lemma3}:}
Recall that $\tilde{\omega}_i=\min\{\omega_i,\sqrt{n}\}$. By definition, the $k$-th moment of $\tilde{\omega}_1$ is equal to
\begin{equation}\nonumber
\begin{aligned}
\mathbb{E}(\tilde{\omega}_1^k)
=& \int_{1}^{+\infty} (\omega_i\wedge\sqrt{n})^k\tau\omega_1^{-\tau-1}d\omega_1\\
=&
\int_{1}^{\sqrt{n}}\tau\omega^{k-\tau-1}d\omega+ \int_{\sqrt{n}}^{+\infty} n^{\frac{k}{2}}\tau\omega^{-\tau-1}d\omega\\
=& 
\frac{\tau}{k-\tau}\omega^{k-\tau}|_1^{\sqrt{n}}+\tau n^{\frac{k}{2}}\frac{1}{(-\tau)}\omega^{-\tau}|_{\sqrt{n}}^{+\infty}\\
=& 
\frac{\tau}{k-\tau}(n^{\frac{k-\tau}{2}}-1)+n^{\frac{k-\tau}{2}}\\
=& 
n^{\frac{k-\tau}{2}}(\frac{\tau}{k-\tau}+1)-\frac{\tau}{k-\tau}\\
=& 
n^{\frac{k-\tau}{2}}\frac{k}{k-\tau}-\frac{\tau}{k-\tau},\ \  k\ne \tau.\\
\end{aligned}
\end{equation}
\qed

\noindent
{\bf Proof of Theorem \ref{theorem:3}}. The proof strategy is similar to the proof of Theorem \ref{theorem:1}. Let  $\mu=\frac{\tau^2}{(\tau-1)^2}$. By Lemma \ref{lemma3},  $\mu_i=\mathbb{E}(d_i)=p\mu$. Simple algebra yields
\begin{eqnarray}\label{rweq1}
\sum_{i=1}^n\left(\frac{d_i}{n}\right)^{2}-\frac{p^{2}\mu^{2}}{n}&=&2\frac{p\mu}{n}\sum_{i=1}^n\frac{d_i-\mu_i}{n}+\sum_{i=1}^n\left(\frac{d_i-\mu_i}{n}\right)^2.
\end{eqnarray}
We now find the order of each term in the right-hand side of (\ref{rweq1}). 

The first term of (\ref{rweq1}) can be decomposed as
\begin{equation}\label{us2}
    \sum_{i=1}^n\frac{d_i-\mu_i}{n}=\frac{\sum_{i\neq j}\left(A_{ij}-p\frac{\tilde{\omega}_i\tilde{\omega}_j}{n}\right)}{n}+\frac{\sum_{i\neq j}\left(p\frac{\tilde{\omega}_i\tilde{\omega}_j}{n}-\frac{p\mu}{n}\right)}{n}.
\end{equation}
Note that $A_{ij}\ (1\leq i<j\leq n)$ are conditionally independent given $W$. Then
\begin{equation}\label{rweq2}
   \mathbb{E}\left[ \frac{\sum_{i< j}\left(A_{ij}-p\frac{\tilde{\omega}_i\tilde{\omega}_j}{n}\right)}{n}\right]^2=\frac{\sum_{i< j}\mathbb{E}\left(A_{ij}-p\frac{\tilde{\omega}_i\tilde{\omega}_j}{n}\right)^2}{n^2}\leq \mathbb{E}\left[p\frac{\tilde{\omega}_i\tilde{\omega}_j}{n}\right]=O\left(\frac{1}{n}\right).
\end{equation}
The second moment of the second term of (\ref{us2}) can be bounded as
\begin{eqnarray}\nonumber
   \mathbb{E}\left[ p\frac{\sum_{i< j}(\tilde{\omega}_i\tilde{\omega}_j-\mu)}{n^2}\right]^2&=&\frac{p^2}{n^4}O\left(\sum_{i\neq j\neq k}\mathbb{E}(\tilde{\omega}_i\tilde{\omega}_j-\mu)(\tilde{\omega}_i\tilde{\omega}_k-\mu)\right)\\ \nonumber
   &&+\frac{p^2}{n^4}O\left(\sum_{i< j}\mathbb{E}(\tilde{\omega}_i\tilde{\omega}_j-\mu)^2\right)\\ \nonumber
   &=&\frac{p^2\mu^2}{n}O\left(n^{\frac{2-\tau}{2}}\right)+\frac{p^2}{n^2}O\left(n^{2-\tau}\right)\\ \label{rweq3}
   &=&O\left(n^{-\frac{\tau}{2}}p^2\mu^2\right),
\end{eqnarray}
where we used Lemma \ref{lemma3} in the second equality. Hence  the first term of (\ref{rweq1}) is $O_P\left(n^{-1-\frac{\tau}{4}}p\mu\right)$.

Now we consider the second  term of (\ref{rweq1}).
By (\ref{rweq2}) and (\ref{rweq3}), we have
\begin{eqnarray}\nonumber
\sum_{i=1}^n\left(\frac{d_i-\mu_i}{n}\right)^2&=&\frac{1}{n^2}\sum_{i=1}^n\left(\sum_{i\neq j}\left(A_{ij}-p\frac{\tilde{\omega}_i\tilde{\omega}_j}{n}\right)\right)^2+\frac{p^2}{n^4}\sum_{i=1}^n\left(\sum_{i\neq j}(\tilde{\omega}_i\tilde{\omega}_j-\mu)\right)^2\\ \nonumber
&&+2\frac{p}{n^3}\sum_{i=1}^n\left(\left[\sum_{i\neq j}\left(A_{ij}-p\frac{\tilde{\omega}_i\tilde{\omega}_j}{n}\right)\right]\left[\sum_{i\neq j}(\tilde{\omega}_i\tilde{\omega}_j-\mu)\right]\right)\\ \label{rweq4}
&=&O_P\left(\frac{1}{n}\right)+\frac{p^2}{n^4}\sum_{i\neq j\neq k}(\tilde{\omega}_i\tilde{\omega}_j-\mu)(\tilde{\omega}_i\tilde{\omega}_k-\mu)+O_P\left(\frac{1}{n^{\frac{1}{2}+\frac{\tau}{4}}}\right),
\end{eqnarray}
where the last term $O_P\left(\frac{1}{n^{\frac{1}{2}+\frac{\tau}{4}}}\right)$ follows from the Cauchy-Schwarz inequality, (\ref{rweq2}) and (\ref{rweq3}). Note that $\mathbb{E}\left[\sum_{i\neq j\neq k}\tilde{\omega}_i^2\tilde{\omega}_j\tilde{\omega}_k\right]\asymp n^{3+\frac{2-\tau}{2}}$. Then
\begin{eqnarray}\nonumber
\sum_{i\neq j\neq k}(\tilde{\omega}_i\tilde{\omega}_j-\mu)(\tilde{\omega}_i\tilde{\omega}_k-\mu)&=&\sum_{i\neq j\neq k}(\tilde{\omega}_i^2\tilde{\omega}_j\tilde{\omega}_k-\tilde{\omega}_i\tilde{\omega}_j\mu-\tilde{\omega}_i\tilde{\omega}_k\mu+\mu^2)\\ \nonumber
&=&(1+o_P(1))\sum_{i\neq j\neq k}\tilde{\omega}_i^2\tilde{\omega}_j\tilde{\omega}_k\\ \nonumber
&=&(1+o_P(1))2\sum_{i<j<k}(\tilde{\omega}_i^2\tilde{\omega}_j\tilde{\omega}_k+\tilde{\omega}_i\tilde{\omega}_j^2\tilde{\omega}_k+\tilde{\omega}_i\tilde{\omega}_j\tilde{\omega}_k^2).
\end{eqnarray}
Hence, by Lemma \ref{Lemma4},
the second term of (\ref{rweq4}) is the leading term and its exact order is $O_P\left(n^{-\frac{\tau}{2}}\right)$. Moreover,
by (\ref{us2}), (\ref{rweq2}) and (\ref{rweq3}), we obtain $\frac{d}{n}=\frac{p\mu}{n}+O_P\left(\frac{p\mu}{n^{1+\frac{\tau}{4}}}\right)$. Then the desired result follows.

\qed

\medskip

\begin{Lemma}\label{Lemma4}
Let $\theta_n=3\mu\left(n^{\frac{2-\tau}{2}}\frac{2}{2-\tau}-\frac{\tau}{2-\tau}\right)$ and
\begin{equation*}\label{rweq6}
U_n=\frac{1}{\binom{n}{3}}\sum_{i<j<k}\left(\tilde{\omega}_i^2\tilde{\omega}_j\tilde{\omega}_k+\tilde{\omega}_i\tilde{\omega}_j^2\tilde{\omega}_k+\tilde{\omega}_i\tilde{\omega}_j\tilde{\omega}_k^2-\theta_n\right).
\end{equation*}
Then 
\begin{equation}\label{ustatistics}
\frac{\sqrt{4-\tau}U_n}{6\mu n^{\frac{1}{2}-\frac{\tau}{4}}}\Rightarrow\mathcal{N}(0,1).
\end{equation}
\end{Lemma}

\medskip

\noindent
{\bf Proof of Lemma \ref{Lemma4}}.   Note that $U_n$ is a U-statistic of order 3. We shall use the asymptotic theory of U-statistics to get the desired result (\ref{ustatistics}). 

Let $\phi(\tilde{\omega}_1,\tilde{\omega}_2,\tilde{\omega}_3)=\tilde{\omega}_1^2\tilde{\omega}_2\tilde{\omega}_3+\tilde{\omega}_1\tilde{\omega}_2^2\tilde{\omega}_3+\tilde{\omega}_1\tilde{\omega}_2\tilde{\omega}_3^2$ and $\phi_1(\tilde{\omega}_1)=\mathbb{E}[\phi(\tilde{\omega}_1,\tilde{\omega}_2,\tilde{\omega}_3)|\tilde{\omega}_1]$. Direct calculation yields $\phi_1(\tilde{\omega}_1)=\tilde{\omega}_1^2\mu+2\tilde{\omega}_1\eta_n$ with $\eta_n=\left(n^{\frac{2-\tau}{2}}\frac{2}{2-\tau}-\frac{\tau}{2-\tau}\right)\frac{\tau}{\tau-1}$.
Note that
\begin{eqnarray*}
\mathbb{E}[U_{n}|\tilde{\omega}_1]
&=\frac{1}{\binom{n}{3}}\sum_{1\leq i<j<k\leq n} \mathbb{E}[\phi(\tilde{\omega}_i,\tilde{\omega}_j,\tilde{\omega}_k)-\theta_n)|\tilde{\omega}_1]=\frac{3}{n}(\phi_1(\tilde{\omega}_1)-\theta_n).
\end{eqnarray*}

Let $\tilde{U}_n=\frac{3}{n}\sum_{i=1}^{n}(\phi_1(\tilde{\omega}_i)-\theta_n)$ and $\sigma_n^2=Var(\tilde{U}_n)$. Then
\begin{eqnarray}\nonumber
\sigma_n^2&=&\frac{9}{n}\mathbb{E}(\phi_1(\tilde{\omega}_1)-\theta_n)^2\\ \nonumber
&=&\frac{9}{n}\Big[\mathbb{E}(\tilde{\omega}_1^4)\mu^2+4\eta_n^2\mathbb{E}(\tilde{\omega}_1^2)+4\mu\eta_n\mathbb{E}(\tilde{\omega}_1^3)-\theta_n^2\Big]\\ \nonumber
&=&(1+o(1))\frac{9}{n}\left[\frac{4\mu^2}{4-\tau}n^{\frac{4-\tau}{2}}+\frac{8\eta_n^2}{2-\tau}n^{\frac{2-\tau}{2}}+\frac{12\mu\eta_n}{3-\tau}n^{\frac{3-\tau}{2}}-\theta_n^2\right]\\
&=&(1+o(1))\frac{36\mu^2}{4-\tau}n^{\frac{4-\tau}{2}-1}.
\end{eqnarray}

Let $Y_i=\frac{3}{n}(\phi_1(\tilde{\omega}_i)-\theta_n)$. Then $Y_i \ (1\leq i\leq n)$ are independent,  $\mathbb{E}(Y_i)=0$ and $\tilde{U}_n=\sum_{i=1}^{n}Y_i$.
Since
\begin{eqnarray}\nonumber
\frac{\sum_{i=1}^n\mathbb{E}(Y_i^4)}{\sigma_n^4}&=&\frac{81}{n^4\sigma^4}\sum_{i=1}^n\mathbb{E}[(\phi_1(\tilde{\omega}_i)-\theta_n)^4]=O\left(\frac{\mathbb{E}(\tilde{\omega}_1^8+\tilde{\omega}_1^4\eta_n^4)}{n^{5-\tau}}\right)\\ \nonumber
&=&O\left(\frac{n^{\frac{8-\tau}{2}}+n^{\frac{4-\tau}{2}+2(2-\tau)})}{n^{5-\tau}}\right)=o(1),
\end{eqnarray}
by the Lyapunov Central Limit Theorem, we get that $\frac{\tilde{U}_n}{\sigma_n}\Rightarrow\mathcal{N}(0,1)$. To finish the proof, it suffices to show $\frac{U_n}{\sigma_n}=\frac{\tilde{U}_n}{\sigma_n}+o_P(1)$. Note that
\begin{equation}\nonumber
\begin{aligned}
\mathbb{E}[\tilde{U}_nU_{n}]
&=\mathbb{E}\left[\frac{3}{n}\sum_{i=1}^{n}(\phi_1(\tilde{\omega}_i)-\theta_n)U_{n}\right]\\
&=\frac{3}{n}\sum_{i=1}^{n}\mathbb{E}[(\phi_1(\tilde{\omega}_i)-\theta_n)\mathbb{E}(U_{n}|\tilde{\omega}_i)]\\
&=\frac{3^2}{n^2}\sum_{i=1}^{n}\mathbb{E}[\phi_1(\tilde{\omega}_i)-\theta_n]^2\\
&=\frac{3^2}{n}\mathbb{E}[\phi_1(\tilde{\omega}_1)-\theta_n]^2\\
&=\frac{3^2}{n}Var(\phi_1(\tilde{\omega}_1))=Var(\tilde{U}_n).
\end{aligned}
\end{equation}
Then
\begin{eqnarray}\nonumber
\mathbb{E}\left[\frac{U_{n}-\tilde{U}_n}{\sigma_n}\right]^2 \nonumber
&=&\frac{1}{\sigma_n^2}\big[\mathbb{E}(U_{n})^2+\mathbb{E}(\tilde{U}_n^2)-2\mathbb{E}(\tilde{U}_nU_{n})\big]\\ \label{ERU}
&=&\frac{1}{\sigma_n^2}[\mathbb{E}(U_{n}^2)-\mathbb{E}(\tilde{U}_n^2)].
\end{eqnarray}

Next, we find $\mathbb{E}(U_{n}^2)$.
\begin{eqnarray}\nonumber
\mathbb{E}(U_{n}^2)&=&\frac{1}{\binom{n}{3}^2}\sum_{\substack{i<j<k,\\  i_1<j_1<k_1}}\mathbb{E}(\phi(\tilde{\omega}_i,\tilde{\omega}_j,\tilde{\omega}_k)-\theta_n)(\phi(\tilde{\omega}_{i_1},\tilde{\omega}_{j_1},\tilde{\omega}_{k_1})-\theta_n)\\ \nonumber
&=&\frac{1}{\binom{n}{3}^2}\sum_{1\leq i<j<k\leq n}\mathbb{E}(\phi(\tilde{\omega}_i,\tilde{\omega}_j,\tilde{\omega}_k)-\theta_n)^2\\ \nonumber
&&+\frac{1}{\binom{n}{3}^2}\sum_{\substack{i<j<k,\\ i_1<j_1<k_1\\
|\{i,j,k\}\cap\{i_1,j_1,k_1\}|=2}}\mathbb{E}(\phi(\tilde{\omega}_i,\tilde{\omega}_j,\tilde{\omega}_k)-\theta_n)(\phi(\tilde{\omega}_{i_1},\tilde{\omega}_{j_1},\tilde{\omega}_{k_1})-\theta_n)\\ \nonumber
&&+\frac{1}{\binom{n}{3}^2}\sum_{\substack{i<j<k,\\ i_1<j_1<k_1\\
|\{i,j,k\}\cap\{i_1,j_1,k_1\}|=1}}\mathbb{E}(\phi(\tilde{\omega}_i,\tilde{\omega}_j,\tilde{\omega}_k)-\theta_n)(\phi(\tilde{\omega}_{i_1},\tilde{\omega}_{j_1},\tilde{\omega}_{k_1})-\theta_n)\\ \label{rweq7}
&=&O\left(\frac{1}{n^{\frac{3}{2}\tau-1}}\right)+O\left(\frac{1}{n^{\tau-1}}\right)+\sigma_n^2(1+o(1)).
\end{eqnarray}
Combining (\ref{ERU}) and (\ref{rweq7}) yields  $\frac{U_n}{\sigma_n}=\frac{\tilde{U}_n}{\sigma_n}+o_p(1)$. Then the proof is complete.

\qed

\noindent
{\bf Proof of Proposition \ref{prop}:}
When $\alpha>1$, the function $f(x)=x^{\alpha}$ is convex for $x>0$. By Jensen inequality, we have
\[\frac{1}{n}\sum_{i=1}^n\left(\frac{d_i}{d}\right)^{\alpha}=\frac{\frac{1}{n}\sum_{i=1}^n\left(\frac{d_i}{n}\right)^{\alpha}}{\left(\frac{1}{n}\sum_{i=1}^n\frac{d_i}{n}\right)^{\alpha}}\geq \frac{\frac{1}{n}\sum_{i=1}^n\left(\frac{d_i}{n}\right)^{\alpha}}{\frac{1}{n}\sum_{i=1}^n\left(\frac{d_i}{n}\right)^{\alpha}}=1.\]
Then $ \mathcal{R}_{\alpha}\in [0,1]$.

When $\alpha\in(0,1)$, the function $f(x)=-x^{\alpha}$ is convex for $x>0$. By Jensen inequality, we have
\[-\frac{1}{n}\sum_{i=1}^n\left(\frac{d_i}{d}\right)^{\alpha}\geq-\left(\frac{1}{n}\sum_{i=1}^n\frac{d_i}{d}\right)^{\alpha}=-1.\]
Then $ \mathcal{R}_{\alpha}\in [0,1]$.

When $\alpha=1$, the function $f(x)=x\log x$ is convex for $x>0$. By Jensen inequality, we have
\[\frac{1}{n}\sum_{i=1}^n
\frac{d_i}{d}\log\left(\frac{d_i}{d}\right)\geq f(1)=0.\]
Then $ \mathcal{R}_{\alpha}\in [0,1]$.

\qed

\section*{Acknowledgement}
The author is grateful to  Editor and anonymous reviewers for valuable comments that significantly improve the manuscript.

\end{document}